\documentclass[11pt,twoside]{article}
\usepackage{amsmath, amsthm, amscd, amsfonts, amssymb, graphicx, color}

\date{}

\begin{document}

\centerline{}

\centerline {\Large{\bf Controlled $g$-atomic subspaces for operators in Hilbert spaces}}

\newcommand{\mvec}[1]{\mbox{\bfseries\itshape #1}}
\centerline{}
\centerline{\textbf{Prasenjit Ghosh}}
\centerline{Department of Pure Mathematics, University of Calcutta,}
\centerline{35, Ballygunge Circular Road, Kolkata, 700019, West Bengal, India}
\centerline{e-mail: prasenjitpuremath@gmail.com}
\centerline{}
\centerline{\textbf{T. K. Samanta}}
\centerline{Department of Mathematics, Uluberia College,}
\centerline{Uluberia, Howrah, 711315,  West Bengal, India}
\centerline{e-mail: mumpu$_{-}$tapas5@yahoo.co.in}

\newtheorem{Theorem}{\quad Theorem}[section]

\newtheorem{definition}[Theorem]{\quad Definition}

\newtheorem{theorem}[Theorem]{\quad Theorem}

\newtheorem{remark}[Theorem]{\quad Remark}

\newtheorem{corollary}[Theorem]{\quad Corollary}

\newtheorem{note}[Theorem]{\quad Note}

\newtheorem{lemma}[Theorem]{\quad Lemma}

\newtheorem{example}[Theorem]{\quad Example}

\newtheorem{result}[Theorem]{\quad Result}
\newtheorem{conclusion}[Theorem]{\quad Conclusion}

\newtheorem{proposition}[Theorem]{\quad Proposition}

\centerline{}

\begin{abstract}
\textbf{\emph{Controlled $g$-atomic subspace for a bounded linear operator is being presented and a characterization has been given.\,We give an example of controlled $K$-$g$-fusion frame.\,We construct a new controlled \,$K$-$g$-fusion frame for the Hilbert space \,$H \,\oplus\, X$\, using the controlled \,$K$-$g$-fusion frames of the Hilbert spaces \,$H$\, and \,$X$.\,Several useful resolutions of the identity operator on a Hilbert space using the theory of controlled $g$-fusion frames have been discussed.\,Frame operator for a pair of controlled $g$-fusion Bessel sequences has been introduced.}}
\end{abstract}



\section{Introduction}
 
\smallskip\hspace{.6 cm} Frame for Hilbert space was first introduced by Duffin and Schaeffer \cite{Duffin} in 1952 to study some fundamental problems in non-harmonic Fourier series.\,Later on, after some decades, frame theory was popularized by Daubechies et al.\,\cite{Daubechies}.For more details on frame theory one can go through the book \cite{O}.\,For special applications, many other types of frames viz.\,$K$-frame \cite{L}, fusion frame \cite{Kutyniok}, \,$g$-frame \cite{Sun}, \,$g$-fusion frame \cite{G, Ahmadi} and \,$K$-$g$-fusion frame \cite{Sadri} etc. were proposed.\,P. Ghosh and T. K. Samanta \cite{P} have studied the stability of dual \,$g$-fusion frame in Hilbert space and also they have discussed generalized atomic subspace for operator in Hilbert space and presented the frame operator for a pair of \,$g$-fusion Bessel sequences \cite{Ghosh}.

Controlled frame is one of the newest generalization of frame.\,I. Bogdanova et al.\,\cite{I} introduced controlled frame for spherical wavelets to get numerically more efficient approximation algorithm.\,Thereafter, P. Balaz \cite{B} developed weighted and controlled frame in Hilbert space.\,In recent times, several generalizations of controlled frame namely, controlled\,$K$-frame \cite{N}, controlled\,$g$-frame \cite{F}, controlled fusion frame \cite{AK}, controlled $g$-fusion frame \cite{HS}, controlled $K$-$g$-fusion frame \cite{GR} etc. have been appeared.\,At present, frame theory has been widely used in signal and image processing, filter bank theory, coding and communications, system modeling and so on.

In this paper,\, we introduce controlled \,$g$-atomic subspace for a bounded linear operator and establish that a family \,$\Lambda_{T\,U}$\, is a controlled \,$g$-atomic subspace if and only if it is a controlled \,$K$-$g$-fusion frame in Hilbert space.\,We construct some useful results about resolution of the identity operator on a Hilbert space using the theory of controlled \,$g$-fusion frame.\,The frame operator for a pair of controlled \,$g$-fusion Bessel sequences are discussed and some of it's properties are going to be established. 

Throughout this paper,\;$H$\; is considered to be a separable Hilbert space with associated inner product \,$\left <\,\cdot \,,\, \cdot\,\right>$\, and \,$\left\{\,H_{j}\,\right\}_{ j \,\in\, J}$\, are the collection of Hilbert spaces, where \,$J$\; is subset of  integers \,$\mathbb{Z}$.\;$I_{H}$\; is the identity operator on \,$H$.\,$\mathcal{B}\,(\,H_{\,1},\, H_{\,2}\,)$\; is a collection of all bounded linear operators from \,$H_{\,1} \,\text{to}\, H_{\,2}$.\,In particular \,$\mathcal{B}\,(\,H\,)$\, denotes the space of all bounded linear operators on \,$H$.\;For \,$S \,\in\, \mathcal{B}\,(\,H\,)$, we denote \,$\mathcal{N}\,(\,S\,)$\; and \,$\mathcal{R}\,(\,S\,)$\, for null space and range of \,$S$, respectively.\,Also, \,$P_{M} \,\in\, \mathcal{B}\,(\,H\,)$\; is the orthonormal projection onto a closed subspace \,$M \,\subset\, H$.\,$\mathcal{G}\,\mathcal{B}\,(\,H\,)$\, denotes the set of all bounded linear operators which have bounded inverse.\,If \,$S,\, R \,\in\, \mathcal{G}\,\mathcal{B}\,(\,H\,)$, then \,$R^{\,\ast},\, R^{\,-\, 1}$\, and \,$S\,R$\, are also belongs to \,$\mathcal{G}\,\mathcal{B}\,(\,H\,)$.\,$\mathcal{G}\,\mathcal{B}^{\,+}\,(\,H\,)$\, is the set of all positive operators in \,$\mathcal{G}\,\mathcal{B}\,(\,H\,)$\, and \,$T,\, U$\, are invertible operators in \,$\mathcal{G}\,\mathcal{B}\,(\,H\,)$.

\section{Preliminaries}
\smallskip\hspace{.6 cm}In this section, we recall some necessary definitions and theorems.

\begin{theorem}(\,Douglas' factorization theorem\,)\,{\cite{Douglas}}\label{th1}
Let \;$S,\, V \,\in\, \mathcal{B}\,(\,H\,)$.\,Then the following conditions are equivalent:
\begin{description}
\item[$(i)$]$\mathcal{R}\,(\,S\,) \,\subseteq\, \mathcal{R}\,(\,V\,)$.
\item[$(ii)$]\;\;$S\, S^{\,\ast} \,\leq\, \lambda^{\,2}\; V\,V^{\,\ast}$\; for some \,$\lambda \,>\, 0$.
\item[$(iii)$]$S \,=\, V\,W$\, for some bounded linear operator \,$W$\, on \,$H$.
\end{description}
\end{theorem}

\begin{theorem}\cite{O}\label{th1.001}
The set \,$\mathcal{S}\,(\,H\,)$\; of all self-adjoint operators on \,$H$\; is a partially ordered set with respect to the partial order \,$\leq$\, which is defined as for \,$R,\,S \,\in\, \mathcal{S}\,(\,H\,)$ 
\[R \,\leq\, S \,\Leftrightarrow\, \left<\,R\,f,\, f\,\right> \,\leq\, \left<\,S\,f,\, f\,\right>\; \;\forall\; f \,\in\, H.\] 
\end{theorem}

\begin{definition}\cite{Kreyzig}
A self-adjoint operator \,$U \,:\, H_{1} \,\to\, H_{1}$\, is called positive if \,$\left<\,U\,x \,,\,  x\,\right> \,\geq\, 0$\, for all \,$x \,\in\, H_{1}$.\;In notation, we can write \,$U \,\geq\, 0$.\;A self-adjoint operator \,$V \,:\, H_{1} \,\to\, H_{1}$\, is called a square root of \,$U$\, if \,$V^{\,2} \,=\, U$.\;If, in addition \,$V \,\geq\, 0$, then \,$V$\, is called positive square root of \,$U$\, and is denoted by \,$V \,=\, U^{1 \,/\, 2}$. 
\end{definition}

\begin{theorem}\cite{Kreyzig}\label{th1.05}
The positive square root \,$V \,:\, H_{1} \,\to\, H_{1}$\, of an arbitrary positive self-adjoint operator \,$U \,:\, H_{1} \,\to\, H_{1}$\, exists and is unique.\;Further, the operator \,$V$\, commutes with every bounded linear operator on \,$H_{1}$\, which commutes with \,$U$.
\end{theorem}

In a complex Hilbert space, every bounded positive operator is self-adjoint and any two bounded positive operators can be commute with each other.

\begin{theorem}\cite{Gavruta}\label{th1.01}
Let \,$M \,\subset\, H$\; be a closed subspace and \,$T \,\in\, \mathcal{B}\,(\,H\,)$.\;Then \,$P_{\,M}\, T^{\,\ast} \,=\, P_{\,M}\,T^{\,\ast}\, P_{\,\overline{T\,M}}$.\;If \,$T$\; is an unitary operator (\,i\,.\,e \,$T^{\,\ast}\, T \,=\, I_{H}$\,), then \,$P_{\,\overline{T\,M}}\;T \,=\, T\,P_{\,M}$.
\end{theorem}

\begin{definition}\cite{Kutyniok}
A family of bounded operators \,$\left\{\,T_{j}\,\right\}_{j \,\in\, J}$\; on \,$H$\; is called a resolution of identity operator on \,$H$\; if for all \,$f \,\in\, H$, we have \,$f \,=\, \sum\limits_{\,j \,\in\, J}\,T_{j}\,(\,f\,)$, provided the series converges unconditionally for all \,$f \,\in\, H$.
\end{definition}

\begin{definition}\cite{Ahmadi}
Let \,$\left\{\,W_{j}\,\right\}_{ j \,\in\, J}$\, be a collection of closed subspaces of \,$H$\; and \,$\left\{\,v_{j}\,\right\}_{ j \,\in\, J}$\; be a collection of positive weights and let \,$\Lambda_{j} \,\in\, \mathcal{B}\,(\,H,\, H_{j}\,)$\; for each \,$j \,\in\, J$.\;Then the family \,$\Lambda \,=\, \{\,\left(\,W_{j},\, \Lambda_{j},\, v_{j}\,\right)\,\}_{j \,\in\, J}$\; is called a generalized fusion frame or a g-fusion frame for \,$H$\; respect to \,$\left\{\,H_{j}\,\right\}_{j \,\in\, J}$\; if there exist constants \,$0 \,<\, A \,\leq\, B \,<\, \infty$\, such that
\begin{equation}\label{eq1}
A \;\left \|\,f \,\right \|^{\,2} \,\leq\, \sum\limits_{\,j \,\in\, J}\,v_{j}^{\,2}\, \left\|\,\Lambda_{j}\,P_{\,W_{j}}\,(\,f\,) \,\right\|^{\,2} \,\leq\, B \; \left\|\, f \, \right\|^{\,2}\; \;\forall\; f \,\in\, H.
\end{equation}
The constants \,$A$\; and \,$B$\; are called the lower and upper bounds of g-fusion frame, respectively.\,If \,$A \,=\, B$\; then \,$\Lambda$\; is called tight g-fusion frame and if \;$A \,=\, B \,=\, 1$\, then we say \,$\Lambda$\; is a Parseval g-fusion frame.\;If  \,$\Lambda$\; satisfies only the right inequality of (\ref{eq1}) it is called a g-fusion Bessel sequence in \,$H$\, with bound \,$B$. 
\end{definition}

Define the space
\[l^{\,2}\left(\,\left\{\,H_{j}\,\right\}_{ j \,\in\, J}\,\right) \,=\, \left \{\,\{\,f_{\,j}\,\}_{j \,\in\, J} \,:\, f_{\,j} \;\in\; H_{j},\; \sum\limits_{\,j \,\in\, J}\, \left \|\,f_{\,j}\,\right \|^{\,2} \,<\, \infty \,\right\}\]
with inner product is given by 
\[\left<\,\{\,f_{\,j}\,\}_{ j \,\in\, J} \,,\, \{\,g_{\,j}\,\}_{ j \,\in\, J}\,\right> \;=\; \sum\limits_{\,j \,\in\, J}\, \left<\,f_{\,j} \,,\, g_{\,j}\,\right>_{H_{j}}.\]\,Clearly \,$l^{\,2}\left(\,\left\{\,H_{j}\,\right\}_{ j \,\in\, J}\,\right)$\; is a Hilbert space with the pointwise operations \cite{Sadri}. 

\begin{definition}\label{def1}\cite{Ghosh}
Let \,$K \,\in\, \mathcal{B}\,(\,H\,)$\; and \,$\left\{\,W_{j}\,\right\}_{ j \,\in\, J}$\; be a collection of closed subspaces of \,$H$, let \,$\left\{\,v_{j}\,\right\}_{ j \,\in\, J}$\, be a collection of positive weights and \,$\Lambda_{j} \,\in\, \mathcal{B}\,(\,H,\, H_{j}\,)$\, for each \,$j \,\in\, J$.\,Then the family \;$\Lambda \,=\, \left\{\,\left(\,W_{j},\, \Lambda_{j},\, v_{j}\,\right)\,\right\}_{j \,\in\, J}$\; is said to be a generalized atomic subspace or g-atomic subspace of \,$H$\; with respect to \,$K$\; if the following statements hold:
\begin{description}
\item[$(i)$]$\Lambda$\; is a g-fusion Bessel sequence in \,$H$.
\item[$(ii)$] For every \,$f \,\in\, H$, there exists \,$\left\{\,f_{\,j}\,\right\}_{j \,\in\, J} \,\in\, l^{\,2}\left(\,\left\{\,H_{j}\,\right\}_{ j \,\in\, J}\,\right)$\; such that 
\[K\,(\,f\,) \,=\, \sum\limits_{\,j \,\in\, J}\, v_{j}\, P_{\,W_{j}}\,\Lambda_{j}^{\,\ast}\, f_{\,j}\; \;\text{and}\; \;\left\|\,\left\{\,f_{\,j}\,\right\}_{j \,\in\, J}\,\right\|_{l^{\,2}\left(\,\left\{\,H_{j}\,\right\}_{ j \,\in\, J}\,\right)} \,\leq\, C\,\|\,f\,\|\] for some \,$C \,>\, 0$. 
\end{description}
\end{definition}

\begin{definition}\cite{Ghosh}
Let \,$\Lambda \,=\, \left\{\,\left(\,V_{i},\, \Lambda_{i},\, v_{i}\,\right)\,\right\}_{i \,\in\, I}$\, and \,$\Lambda^{\,\prime} \,=\, \left\{\,\left(\,V^{\,\prime}_{i},\, \Lambda^{\,\prime}_{i},\, v^{\,\prime}_{i}\,\right)\,\right\}_{i \,\in\, I}$\, be two g-fusion Bessel sequences in \,$H$\; with bounds \,$D_{\,1}$\, and \,$D_{\,2}$, respectively.\;Then the operator \,$S_{\Lambda\,\Lambda^{\,\prime}} \,:\, H \,\to\, H$, defined by
\[S_{\Lambda\,\Lambda^{\,\prime}}\,(\,f\,) \;=\; \sum\limits_{\,i \,\in\, I}\,v_{\,i}\,v^{\,\prime}_{\,i}\,P_{\,V_{i}}\,\Lambda_{i}^{\,\ast}\,\Lambda_{i}^{\,\prime} \,P_{\,V^{\,\prime}_{i}}\,(\,f\,)\;\; \;\forall\; f \,\in\, H,\] is called the frame operator for the pair of g-fusion Bessel sequences \,$\Lambda$\, and \,$\Lambda^{\,\prime}$.
\end{definition}

\begin{definition}\cite{HS}
Let \,$\left\{\,W_{j}\,\right\}_{ j \,\in\, J}$\, be a collection of closed subspaces of \,$H$\, and \,$\left\{\,v_{j}\,\right\}_{ j \,\in\, J}$\, be a collection of positive weights.\,Let \,$\left\{\,H_{j}\,\right\}_{ j \,\in\, J}$\, be a sequence of Hilbert spaces, \,$T,\, U \,\in\, \mathcal{G}\,\mathcal{B}\,(\,H\,)$\, and \,$\Lambda_{j} \,\in\, \mathcal{B}\,(\,H,\, H_{j}\,)$\, for each \,$j \,\in\, J$.\,Then the family \,$\Lambda_{T\,U} \,=\, \left\{\,\left(\,W_{j},\, \Lambda_{j},\, v_{j}\,\right)\,\right\}_{j \,\in\, J}$\, is a \,$(\,T,\,U\,)$-controlled $g$-fusion frame for \,$H$\, if there exist constants \,$0 \,<\, A \,\leq\, B \,<\, \infty$\, such that 
\begin{equation}\label{eqn1.1}
A\,\|\,f\,\|^{\,2} \,\leq\, \sum\limits_{\,j \,\in\, J}\, v^{\,2}_{j}\,\left<\,\Lambda_{j}\,P_{\,W_{j}}\,U\,f,\,  \Lambda_{j}\,P_{\,W_{j}}\,T\,f\,\right> \,\leq\, \,B\,\|\,f \,\|^{\,2}\; \;\forall\; f \,\in\, H.
\end{equation}
If \,$A \,=\, B$\, then \,$\Lambda_{T\,U}$\, is called \,$(\,T,\,U\,)$-controlled tight g-fusion frame and if \,$A \,=\, B \,=\, 1$\, then we say \,$\Lambda_{T\,U}$\, is a \,$(\,T,\,U\,)$-controlled Parseval g-fusion frame.\,If \,$\Lambda_{T\,U}$\, satisfies only the right inequality of (\ref{eqn1.1}) it is called a \,$(\,T,\,U\,)$-controlled g-fusion Bessel sequence in \,$H$.    
\end{definition}

\begin{definition}\cite{HS}
Let \,$\Lambda_{T\,U}$\, be a \,$(\,T,\,U\,)$-controlled g-fusion Bessel sequence in \,$H$\, with a bound \,$B$.\,The synthesis operator \,$T_{C} \,:\, \mathcal{K}_{\,\Lambda_{j}} \,\to\, H$\, is defined as 
\[T_{C}\,\left(\,\left\{\,v_{\,j}\,\left(\,T^{\,\ast}\,P_{\,W_{j}}\, \Lambda_{j}^{\,\ast}\,\Lambda_{j}\,P_{\,W_{j}}\,U\,\right)^{1 \,/\, 2}\,f\,\right\}_{j \,\in\, J}\,\right) \,=\, \sum\limits_{\,j \,\in\, J}\,v^{\,2}_{j}\,T^{\,\ast}\,P_{\,W_{j}}\, \Lambda_{j}^{\,\ast}\,\Lambda_{j}\,P_{\,W_{j}}\,U\,f,\]for all \,$f \,\in\, H$\, and the analysis operator \,$T^{\,\ast}_{C} \,:\, H \,\to\, \mathcal{K}_{\,\Lambda_{j}}$\,is given by 
\[T_{C}^{\,\ast}\,f \,=\,  \left\{\,v_{\,j}\,\left(\,T^{\,\ast}\,P_{\,W_{j}}\, \Lambda_{j}^{\,\ast}\,\Lambda_{j}\,P_{\,W_{j}}\,U\,\right)^{1 \,/\, 2}\,f\,\right\}_{j \,\in\, J}\; \;\forall\; f \,\in\, H,\]
where 
\[\mathcal{K}_{\,\Lambda_{j}} \,=\, \left\{\,\left\{\,v_{\,j}\,\left(\,T^{\,\ast}\,P_{\,W_{j}}\, \Lambda_{j}^{\,\ast}\,\Lambda_{j}\,P_{\,W_{j}}\,U\,\right)^{1 \,/\, 2}\,f\,\right\}_{j \,\in\, J} \,:\, f \,\in\, H\,\right\} \,\subset\, l^{\,2}\left(\,\left\{\,H_{j}\,\right\}_{ j \,\in\, J}\,\right).\]
The frame operator \,$S_{C} \,:\, H \,\to\, H$\; is defined as follows:
\[S_{C}\,f \,=\, T_{C}\,T_{C}^{\,\ast}\,f \,=\, \sum\limits_{\,j \,\in\, J}\, v_{j}^{\,2}\,T^{\,\ast}\,P_{\,W_{j}}\, \Lambda_{j}^{\,\ast}\,\Lambda_{j}\,P_{\,W_{j}}\,U\,f\; \;\forall\; f \,\in\, H\]and it is easy to verify that 
\[\left<\,S_{C}\,f,\, f\,\right> \,=\, \sum\limits_{\,j \,\in\, J}\, v^{\,2}_{j}\,\left<\,\Lambda_{j}\,P_{\,W_{j}}\,U\,f,\,  \Lambda_{j}\,P_{\,W_{j}}\,T\,f\,\right>\; \;\forall\; f \,\in\, H.\]
Furthermore, if \,$\Lambda_{T\,U}$\, is a \,$(\,T,\,U\,)$-controlled g-fusion frame with bounds \,$A$\, and \,$B$\, then \,$A\,I_{\,H} \,\leq\,S_{C} \,\leq\, B\,I_{H}$.\,Hence, \,$S_{C}$\, is bounded, invertible, self-adjoint and positive linear operator.\,It is easy to verify that \,$B^{\,-1}\,I_{H} \,\leq\, S_{C}^{\,-1} \,\leq\, A^{\,-1}\,I_{H}$.
\end{definition}

\begin{definition}\cite{GR}
Let \,$K \,\in\, \mathcal{B}\,(\,H\,)$\, and \,$\left\{\,W_{j}\,\right\}_{ j \,\in\, J}$\, be a collection of closed subspaces of \,$H$\, and \,$\left\{\,v_{j}\,\right\}_{ j \,\in\, J}$\, be a collection of positive weights.\,Let \,$\left\{\,H_{j}\,\right\}_{ j \,\in\, J}$\, be a sequence of Hilbert spaces, \,$T,\, U \,\in\, \mathcal{G}\,\mathcal{B}\,(\,H\,)$\, and \,$\Lambda_{j} \,\in\, \mathcal{B}\,(\,H,\, H_{j}\,)$\, for each \,$j \,\in\, J$.\,Then the family \,$\Lambda_{T\,U} \,=\, \left\{\,\left(\,W_{j},\, \Lambda_{j},\, v_{j}\,\right)\,\right\}_{j \,\in\, J}$\, is a \,$(\,T,\,U\,)$-controlled $K$-$g$-fusion frame for \,$H$\, if there exist constants \,$0 \,<\, A \,\leq\, B \,<\, \infty$\, such that 
\[A\,\|\,K^{\,\ast}\,f\,\|^{\,2} \,\leq\, \sum\limits_{\,j \,\in\, J}\, v^{\,2}_{j}\,\left<\,\Lambda_{j}\,P_{\,W_{j}}\,U\,f,\,  \Lambda_{j}\,P_{\,W_{j}}\,T\,f\,\right> \,\leq\, \,B\,\|\,f \,\|^{\,2}\; \;\forall\; f \,\in\, H.\]
\end{definition}

Throughout this paper, \,$\Lambda_{T\,U}$\, denotes the family \,$\left\{\,\left(\,W_{j},\, \Lambda_{j},\, v_{j}\,\right)\,\right\}_{j \,\in\, J}$.

\section{Controlled $g$-atomic subspace}

\smallskip\hspace{.6 cm} In this section, we define controlled generalized atomic subspace or controlled \,$g$-atomic subspace of a Hilbert space with respect to a bounded linear operator.\,At the end of this section, a controlled \,$K$-$g$-fusion frame for the Hilbert space \,$H \,\oplus\, X$, using the controlled \,$K$-$g$-fusion frames of the Hilbert spaces \,$H$\, and \,$X$\, will be constructed. 

\begin{definition}\label{def1}
Let \,$K \,\in\, \mathcal{B}\,(\,H\,)$\, and \,$\left\{\,W_{j}\,\right\}_{ j \,\in\, J}$\, be a family of closed subspaces of \,$H$\, and \,$\left\{\,v_{j}\,\right\}_{ j \,\in\, J}$\, be a family of positive weights.\,Let \,$\left\{\,H_{j}\,\right\}_{ j \,\in\, J}$\, be a sequence of Hilbert spaces, \,$T,\, U \,\in\, \mathcal{G}\,\mathcal{B}\,(\,H\,)$\, and \,$\Lambda_{j} \,\in\, \mathcal{B}\,(\,H,\, H_{j}\,)$\, for each \,$j \,\in\, J$.\,Then the family \,$\Lambda_{T\,U} \,=\, \left\{\,\left(\,W_{j},\, \Lambda_{j},\, v_{j}\,\right)\,\right\}_{j \,\in\, J}$\, is said to be a generalized atomic subspace controlled by the operators \,$T,\, U$\, or \,$(\,T,\,U\,)$-controlled g-atomic subspace of \,$H$\, with respect to \,$K$\, if the following statements hold:
\begin{description}
\item[$(i)$]\,$\Lambda_{T\,U}$\, is a \,$(\,T,\,U\,)$-controlled g-fusion Bessel sequence in \,$H$.
\item[$(ii)$]\, For every \,$f \,\in\, H$, there exists \,$\left\{\,v_{\,j}\,\left(\,T^{\,\ast}\,P_{\,W_{j}}\, \Lambda_{j}^{\,\ast}\; \Lambda_{j}\, P_{\,W_{j}}\,U\,\right)^{1 \,/\, 2}\,f\,\right\}_{j \,\in\, J} \,\in\, \mathcal{K}_{\Lambda_{j}}^{\,2}$\, such that \,$K\,f \,=\, \sum\limits_{\,j \,\in\, J}\, v^{\,2}_{j}\,T^{\,\ast}\,P_{\,W_{j}}\, \Lambda_{j}^{\,\ast}\, \Lambda_{j}\, P_{\,W_{j}}\,U\,f$, and
\[\left\|\,\left\{\,v_{\,j}\,\left(\,T^{\,\ast}\,P_{\,W_{j}}\, \Lambda_{j}^{\,\ast}\; \Lambda_{j}\, P_{\,W_{j}}\,U\,\right)^{1 \,/\, 2}\,f\,\right\}_{j \,\in\, J} \,\right\|_{\mathcal{K}_{\Lambda_{j}}^{\,2}} \; \leq\; C\; \|\,f\,\|\]
for some \,$C \,>\, 0$. 
\end{description}
\end{definition}

The next theorem provides a characterization of a controlled \,$g$-atomic subspace. 

\begin{theorem}\label{thm1}
Let \,$K \,\in\, \mathcal{B}\,(\,H\,)$\, and \,$\left\{\,W_{j}\,\right\}_{ j \,\in\, J}$\, be a family of closed subspaces of \,$H$, let \,$\left\{\,v_{j}\,\right\}_{ j \,\in\, J}$\, be a family of positive weights.\,Let \,$T,\, U \,\in\, \mathcal{G}\,\mathcal{B}\,(\,H\,)$\, and \,$\Lambda_{j} \,\in\, \mathcal{B}\,(\,H \,,\, H_{j}\,)$\; for each \,$j \,\in\, J$.\;Then the following statements are equivalent:
\begin{description}
\item[$(i)$]\,\,$\Lambda_{T\,U} \,=\, \left\{\,\left(\,W_{j},\, \Lambda_{j},\, v_{j}\,\right)\,\right\}_{j \,\in\, J}$\, is a \,$(\,T,\,U\,)$-controlled g-atomic subspace of \,$H$\, with respect to \,$K$.
\item[$(ii)$] \,$\Lambda_{T\,U}$\, is a \,$(\,T,\,U\,)$-controlled K-g-fusion frame for \,$H$.
\end{description}
\end{theorem}

\begin{proof}$(i) \,\Rightarrow\, (ii)$
Suppose \,$\Lambda_{T\,U}$\, is a \,$(\,T,\,U\,)$-controlled $g$-atomic subspace of \,$H$\, with respect to \,$K$.\,Then \,$\Lambda_{T\,U}$\, is a \,$(\,T,\,U\,)$-controlled \,$g$-fusion Bessel sequence in \,$H$.\,Then for each \,$f \,\in\, H$, there exists constant \,$B \,>\, 0$\, such that
\[\sum\limits_{\,j \,\in\, J}\, v^{\,2}_{j}\,\left<\, \Lambda_{j}\, P_{\,W_{j}}\,U\,f,\,  \Lambda_{j}\, P_{\,W_{j}}\,T\,f\,\right> \,\leq\, B\,\|\,f\,\|^{\,2}.\] 
Now, for any \,$f \,\in\, H$, we have
\[\left \|\,K^{\,\ast}\,f \,\right \| \,=\, \sup\limits_{\,\|\,g\,\| \,=\, 1}\,\left|\,\left<\,K^{\,\ast}\,f \,,\, g\,\right>\,\right| \,=\, \sup\limits_{\,\|\,g\,\| \,=\, 1}\,\left|\,\left<\,f \,,\, K\,g\,\right>\,\right|,\]
by definition (\,\ref{def1}\,), for \,$g \,\in\, H$, there exists \,$\left\{\,v_{\,j}\,\left(\,T^{\,\ast}\,P_{\,W_{j}}\, \Lambda_{j}^{\,\ast}\,\Lambda_{j}\,P_{\,W_{j}}\,U\,\right)^{1 \,/\, 2}\,g\,\right\}_{j \,\in\, J} \,\in\, \mathcal{K}_{\Lambda_{j}}^{\,2}$\, such that \,$K\,g \,=\, \sum\limits_{\,j \,\in\, J}\, v^{\,2}_{j}\,T^{\,\ast}\,P_{\,W_{j}}\, \Lambda_{j}^{\,\ast}\, \Lambda_{j}\, P_{\,W_{j}}\,U\,g$, and 
\[\left\|\,\left\{\,v_{\,j}\,\left(\,T^{\,\ast}\,P_{\,W_{j}}\, \Lambda_{j}^{\,\ast}\, \Lambda_{j}\, P_{\,W_{j}}\,U\,\right)^{1 \,/\, 2}\,g\,\right\}_{j \,\in\, J} \,\right\|_{\mathcal{K}_{\Lambda_{j}}^{\,2}} \,\leq\, C\, \|\,g\,\|\] 
for some \,$C \,>\, 0$.\,Thus
\begin{align*}
&\left \|\,K^{\,\ast}\,f \,\right \|^{\,2} \,=\, \sup\limits_{\,\|\,g\,\| \,=\, 1}\,\left|\,\left<\,f,\, \sum\limits_{\,j \,\in\, J}\, v^{\,2}_{j}\,T^{\,\ast}\,P_{\,W_{j}}\, \Lambda_{j}^{\,\ast}\, \Lambda_{j}\, P_{\,W_{j}}\,U\,g\,\right>\,\right|^{\,2}\\
&\,=\, \sup\limits_{\,\|\,g\,\| \,=\, 1}\,\left|\,\sum\limits_{\,j \,\in\, J}\,\left<\,v_{\,j}\,\left(\,T^{\,\ast}\,P_{\,W_{j}}\, \Lambda_{j}^{\,\ast}\, \Lambda_{j}\, P_{\,W_{j}}\,U\,\right)^{1 \,/\, 2}\,f,\, v_{\,j}\,\left(\,T^{\,\ast}\,P_{\,W_{j}}\, \Lambda_{j}^{\,\ast}\, \Lambda_{j}\, P_{\,W_{j}}\,U\,\right)^{1 \,/\, 2}\,g\,\right>\,\right|^{\,2}\\
&\leq\, \sup\limits_{\,\|\,g\,\| \,=\, 1}\,\sum\limits_{\,j \,\in\, J}\,\left\|\,v_{\,j}\,\left(\,T^{\,\ast}\,P_{\,W_{j}}\, \Lambda_{j}^{\,\ast}\, \Lambda_{j}\, P_{\,W_{j}}\,U\,\right)^{1 \,/\, 2}\,f\,\right\|^{\,2}\,\left\|\,\left\{\,v_{\,j}\,\left(\,T^{\,\ast}\,P_{\,W_{j}}\, \Lambda_{j}^{\,\ast}\, \Lambda_{j}\, P_{\,W_{j}}\,U\,\right)^{1 \,/\, 2}\,g\,\right\}_{j \,\in\, J} \,\right\|^{\,2}_{\mathcal{K}_{\Lambda_{j}}^{\,2}}\\
&\leq\,\sup\limits_{\,\|\,g\,\| \,=\, 1}\,\sum\limits_{\,j \,\in\, J}\, v^{\,2}_{j}\,\left<\, \Lambda_{j}\, P_{\,W_{j}}\,U\,f,\,  \Lambda_{j}\, P_{\,W_{j}}\,T\,f\,\right>\;C^{\,2}\,\|\,g\,\|^{\,2}\\
&\Rightarrow\, \dfrac{1}{C^{\,2}}\,\left \|\,K^{\,\ast}\,f \,\right \|^{\,2} \,\leq\, \sum\limits_{\,j \,\in\, J}\, v^{\,2}_{j}\,\left<\, \Lambda_{j}\, P_{\,W_{j}}\,U\,f,\,  \Lambda_{j}\, P_{\,W_{j}}\,T\,f\,\right>.
\end{align*}
Thus, \,$\Lambda_{T\,U}$\, is a \,$(\,T,\,U\,)$-controlled $K$-$g$-fusion frame for \,$H$\, with bounds \,$\dfrac{1}{C^{\,2}}$\, and \,$B$.\\\\
$(ii) \,\Rightarrow\, (i)$\, Suppose that \,$\Lambda_{T\,U}$\, is a \,$(\,T,\,U\,)$-controlled $K$-$g$-fusion frame for \,$H$\, with the corresponding synthesis operator \,$T_{C}$.\,Then obviously \,$\Lambda_{T\,U}$\, is a \,$(\,T,\,U\,)$-controlled \,$g$-fusion Bessel sequence in \,$H$.\,Now, for each \,$f \,\in\, H$,
\[A\, \left \|\,K^{\,\ast}\,f \,\right \|^{\,2} \,\leq\, \sum\limits_{\,j \,\in\, J}\, v^{\,2}_{j}\,\left<\, \Lambda_{j}\, P_{\,W_{j}}\,U\,f,\,  \Lambda_{j}\, P_{\,W_{j}}\,T\,f\,\right> \,=\, \|\,T_{C}^{\,\ast}\,f\,\|^{\,2}\]
gives \,$A \,K\,K^{\,\ast} \,\leq\, T_{C}\,T_{C}^{\,\ast}$\, and by Theorem (\ref{th1}), there exists \,$L \,\in\, \mathcal{B}\,\left(\,H,\,\mathcal{K}_{\Lambda_{j}}^{\,2}\,\right)$\, such that \,$K \,=\, T_{C}\,L$.\,Define 
\[L\,(\,f\,) \,=\, \left\{\,v_{\,j}\,\left(\,T^{\,\ast}\,P_{\,W_{j}}\, \Lambda_{j}^{\,\ast}\, \Lambda_{j}\, P_{\,W_{j}}\,U\,\right)^{1 \,/\, 2}\,f\,\right\}_{j \,\in\, J},\] for every \,$f \,\in\, H$.\,Then, for each \,$f \,\in\, H$, we have
\begin{align*}
&K \,(\,f\,) \,=\, T_{C}\,L\,(\,f\,) \,=\, \sum\limits_{\,j \,\in\, J}\, v^{\,2}_{j}\,T^{\,\ast}\,P_{\,W_{j}}\, \Lambda_{j}^{\,\ast}\, \Lambda_{j}\, P_{\,W_{j}}\,U\,f,\; \;\text{and}\\ 
&\left\|\,\left\{\,v_{\,j}\,\left(\,T^{\,\ast}\,P_{\,W_{j}}\, \Lambda_{j}^{\,\ast}\, \Lambda_{j}\, P_{\,W_{j}}\,U\,\right)^{1 \,/\, 2}\,f\,\right\}_{j \,\in\, J} \,\right\|_{\mathcal{K}_{\Lambda_{\,j}}^{\,2}} \,=\, \left\|\,L\,(\,f\,)\,\right\|_{\mathcal{K}_{\Lambda_{j}}^{\,2}} \,\leq\, C\; \|\,f\,\|,
\end{align*} 
where \,$C \,=\, \|\,L\,\|$.\,Thus, \,$\Lambda_{T\,U}$\, is a \,$(\,T,\,U\,)$-controlled $g$-atomic subspace of \,$H$\, with respect to \,$K$.
\end{proof}

\begin{remark}
Let \,$K_{1},\, K_{2} \,\in\, \mathcal{B}\,(\,H\,)$.\,If \,$\Lambda_{T\,U}$\, is an \,$(\,T,\,U\,)$-controlled $g$-atomic subspace of \,$H$\, with respect to \,$K_{1}$\, and \,$K_{2}$, then for any scalars \,$\alpha$\, and \,$\beta$, \,$\Lambda_{T\,U}$\, is an \,$(\,T,\,U\,)$-controlled $g$-atomic subspace of \,$H$\, with respect to \,$\alpha\,K_{1} \,+\, \beta\,K_{2}$\, and \,$K_{1}\,K_{2}$.  
\end{remark}

In the next theorem, we will see that every controlled \,$g$-fusion frame is a controlled \,$g$-atomic subspace with respect to it's frame operator.

\begin{theorem}
Let \,$\Lambda_{T\,U} \,=\, \left\{\,\left(\,W_{j},\, \Lambda_{j},\, v_{j}\,\right)\,\right\}_{j \,\in\, J}$\, be a \,$(\,T,\,U\,)$-controlled g-fusion frame for \,$H$.\,Then \,$\Lambda_{T\,U}$\, is a \,$(\,T,\,U\,)$-controlled $g$-atomic subspace of \,$H$\, with respect to it's frame operator \,$S_{C}$.
\end{theorem}

\begin{proof}
Since \,$\Lambda_{T\,U}$\, is a \,$(\,T,\,U\,)$-controlled $g$-fusion frame for \,$H$, we have \,$\mathcal{R}\,\left(T_{C}\,\right) \,=\, H \,=\, \mathcal{R}\,\left(\,S_{C}\,\right)$, so by Theorem (\ref{th1}), there exists some \,$\alpha \,>\, 0$\; such that \,$\alpha \,S_{C}\, S_{C}^{\,\ast} \,\leq\, T_{C}\,T_{C}^{\,\ast}$\; and therefore for each \,$f \,\in\, H$, we have
\[\alpha\,\left \|\,S_{C}^{\,\ast}\,f \,\right \|^{\,2} \,\leq\, \|\,T_{C}^{\,\ast}\,f\,\|^{\,2} \,=\, \sum\limits_{\,j \,\in\, J}\, v^{\,2}_{j}\,\left<\, \Lambda_{j}\, P_{\,W_{j}}\,U\,f,\,  \Lambda_{j}\, P_{\,W_{j}}\,T\,f\,\right>  \,\leq\, B\,\left\|\, f \, \right\|^{\,2}.\] Thus, \,$\Lambda_{T\,U}$\, is a \,$(\,T,\,U\,)$-controlled $S_{C}$-$g$-fusion frame for \,$H$\, and hence by Theorem (\ref{thm1}), \,$\Lambda_{T\,U}$\, is a \,$(\,T,\,U\,)$-controlled $g$-atomic subspace of \,$H$\, with respect to \,$S_{C}$.
\end{proof}

Now, we construct another controlled \,$g$-atomic subspace from given controlled \,$g$-atomic subspaces under some sufficient conditions. 

\begin{theorem}\label{thm2}
Let \,$\Lambda_{T\,U} \,=\, \left\{\,\left(\,W_{j},\, \Lambda_{j},\, v_{j}\,\right)\,\right\}_{j \,\in\, J}$\, and \,$\Gamma_{T\,U} \,=\, \left\{\,\left(\,W_{j},\, \Gamma_{j},\, v_{j}\,\right)\,\right\}_{j \,\in\, J}$\, be two \,$(\,T,\,U\,)$-controlled $g$-atomic subspaces of \,$H$\, with respect to \,$K \,\in\, \mathcal{B}\,(\,H\,)$\, and let \,$V,\, W \,\in\, \mathcal{B}\,(\,H\,)$\, such that \,$V \,+\, W$\, be an invertible operator on \,$H$\, with \,$K\, \left(\,V \,+\, W \,\right) \,=\, \left(\,V \,+\, W \,\right)\, K$.\,Suppose \,$(\,V \,+\, W\,)^{\,\ast}$\, commutes with \,$T$\, and \,$U$.\,Then
\[\Delta_{T\,U} \,=\, \left\{\,\left(\,(\,V \,+\, W\,)\,W_{j},\; \left(\,\Lambda_{j} \,+\, \Gamma_{j}\,\right)\,P_{\,W_{j}}\,\left(\,V \,+\, W \,\right)^{\,\ast},\; v_{j}\,\right)\,\right\}_{j \,\in\, J}\] is a \,$(\,T,\,U\,)$-controlled $g$-atomic subspace of \,$H$\, with respect to \,$K$, if for all \,$j$,
\[\left<\,\Gamma_{j}\,P_{\,W_{j}}\,(\,V \,+\, W\,)^{\,\ast}\,U\,f,\, \Lambda_{j}\,P_{\,W_{j}}\,(\,V \,+\, W\,)^{\,\ast}\,T\,f\,\,\right> \,=\, 0,\] 
\begin{equation}\label{eq2.1}
 \left<\,\Lambda_{j}\,P_{\,W_{j}}\,(\,V \,+\, W\,)^{\,\ast}\,U\,f,\, \Gamma_{j}\,P_{\,W_{j}}\,(\,V \,+\, W\,)^{\,\ast}\,T\,f\,\right> \,=\, 0,\; f \,\in\, H.
\end{equation}
\end{theorem}

\begin{proof}
Since \,$\Lambda_{T\,U}$\, and \,$\Gamma_{T\,U}$\, are \,$(\,T,\,U\,)$-controlled\,$g$-atomic subspaces with respect to \,$K$, by Theorem (\ref{thm1}), they are \,$(\,T,\,U\,)$-controlled \,$K$-$g$-fusion frames for \,$H$.\,So, for each \,$f \,\in\, H$, there exist positive constants \,$(\,A,\, B\,)$\, and \,$(\,C,\, D\,)$\, such that
\begin{align*}
&A\,\|\,K^{\,\ast}\,f\,\|^{\,2} \,\leq\, \sum\limits_{\,j \,\in\, J}\, v^{\,2}_{j}\,\left<\,\Lambda_{j}\,P_{\,W_{j}}\,U\,f,\,  \Lambda_{j}\,P_{\,W_{j}}\,T\,f\,\right> \,\leq\, \,B\,\|\,f \,\|^{\,2},\;\text{and}\\
&C\,\|\,K^{\,\ast}\,f\,\|^{\,2} \,\leq\, \sum\limits_{\,j \,\in\, J}\, v^{\,2}_{j}\,\left<\,\Gamma_{j}\,P_{\,W_{j}}\,U\,f,\,  \Gamma_{j}\,P_{\,W_{j}}\,T\,f\,\right> \,\leq\, D\,\|\,f \,\|^{\,2}. 
\end{align*} 
Also, \,$V \,+\, W$\, is invertible, so
\[\left\|\,K^{\,\ast}\,f\,\right\|^{\,2} \,=\, \left\|\,\left(\,\left(\,V \,+\, W\,\right)^{\,-\, 1}\,\right)^{\,\ast}\,(\,V \,+\, W\,)^{\,\ast}\,K^{\,\ast}\,f\,\right\|^{\,2}\]
\begin{equation}\label{eqn4}
\hspace{1.9cm}\leq\, \left\|\,\left(\,V \,+\, W\right)^{\,-\, 1}\,\right\|^{\,2}\, \left\|\,(\,V \,+\, W\,)^{\,\ast}\,K^{\,\ast}\,f\,\right\|^{\,2}.
\end{equation} 
Take \,$R \,=\, V \,+\, W$.\,Now, for each \,$f \,\in\, H$, using Theorem \ref{th1.01}, we have
\begin{align*}
&\sum\limits_{\,j \,\in\, J}\,v_{j}^{\,2}\,\left<\,\left(\,\Lambda_{j} \,+\, \Gamma_{j}\,\right)\,P_{\,W_{j}}\,R^{\,\ast}\, P_{\,R\,W_{j}}\,U\,f,\, \left(\,\Lambda_{j} \,+\, \Gamma_{j}\,\right)\,P_{\,W_{j}}\,R^{\,\ast}\, P_{\,R\,W_{j}}\,T\,f\,\right>\\  
&=\, \sum\limits_{\,j \,\in\, J}\,v_{j}^{\,2}\,\left<\,\left(\,\Lambda_{j} \,+\, \Gamma_{j}\,\right)\,P_{\,W_{j}}\,R^{\,\ast}\,U\,f,\, \left(\,\Lambda_{j} \,+\, \Gamma_{j}\,\right)\,P_{\,W_{j}}\,R^{\,\ast}\,T\,f\,\right>\\
&=\, \sum\limits_{\,j \,\in\, J}\,v_{j}^{\,2}\,\left<\,\left(\,\Lambda_{j} \,+\, \Gamma_{j}\,\right)\,P_{\,W_{j}}\,U\,R^{\,\ast}\,f,\, \left(\,\Lambda_{j} \,+\, \Gamma_{j}\,\right)\,P_{\,W_{j}}\,T\,R^{\,\ast}\,f\,\right>\\
&=\, \sum\limits_{\,j \,\in\, J}\,v_{j}^{\,2}\,\left\{\,\left<\,\Lambda_{j}\,P_{\,W_{j}}\,U\,R^{\,\ast}\,f,\, \Lambda_{j}\,P_{\,W_{j}}\,T\,R^{\,\ast}\,f\,\right> \,+\,\left<\,\Gamma_{j}\,P_{\,W_{j}}\,U\,R^{\,\ast}\,f,\, \Gamma_{j}\,P_{\,W_{j}}\,T\,R^{\,\ast}\,f\,\right>\,\right\}\,+ \\
&+\,\sum\limits_{\,j \,\in\, J}\,v_{j}^{\,2}\,\left\{\,\left<\,\Lambda_{j}\,P_{\,W_{j}}\,U\,R^{\,\ast}\,f,\, \Gamma_{j}\,P_{\,W_{j}}\,T\,R^{\,\ast}\,f\,\right> \,+\,\left<\,\Gamma_{j}\,P_{\,W_{j}}\,U\,R^{\,\ast}\,f,\, \Lambda_{j}\,P_{\,W_{j}}\,T\,R^{\,\ast}\,f\,\right>\,\right\}\\
&=\, \sum\limits_{\,j \,\in\, J}\,v_{j}^{\,2}\,\left\{\,\left<\,\Lambda_{j}\,P_{\,W_{j}}\,U\,R^{\,\ast}\,f,\, \Lambda_{j}\,P_{\,W_{j}}\,T\,R^{\,\ast}\,f\,\right> \,+\,\left<\,\Gamma_{j}\,P_{\,W_{j}}\,U\,R^{\,\ast}\,f,\, \Gamma_{j}\,P_{\,W_{j}}\,T\,R^{\,\ast}\,f\,\right>\,\right\}\\
&\leq\, B\, \left\|\,R^{\,\ast}\,f\,\right\|^{\,2} \,+\, D\,\left\|\,R^{\,\ast}\,f\,\right\|^{\,2} \,=\, \left(\,B \,+\, D\,\right)\,\left\|\,\left(\,V \,+\, W\,\right)^{\,\ast}\,f\,\right\|^{\,2}\\
&\leq\, \left(\,B \,+\, D\,\right)\, \left\|\,V \,+\, W\,\right\|^{\,2}\, \|\,f\,\|^{\,2}\; \;[\;\text{as}\; \,V \,+\, W \;\text{is bounded}\;].
\end{align*}
On the other hand,
\begin{align*}
&\sum\limits_{\,j \,\in\, J}\,v_{j}^{\,2}\,\left<\,\left(\,\Lambda_{j} \,+\, \Gamma_{j}\,\right)\,P_{\,W_{j}}\,R^{\,\ast}\, P_{\,R\,W_{j}}\,U\,f,\, \left(\,\Lambda_{j} \,+\, \Gamma_{j}\,\right)\,P_{\,W_{j}}\,R^{\,\ast}\, P_{\,R\,W_{j}}\,T\,f\,\right>\\
&=\, \sum\limits_{\,j \,\in\, J}\,v_{j}^{\,2}\,\left\{\,\left<\,\Lambda_{j}\,P_{\,W_{j}}\,U\,R^{\,\ast}\,f,\, \Lambda_{j}\,P_{\,W_{j}}\,T\,R^{\,\ast}\,f\,\right> \,+\,\left<\,\Gamma_{j}\,P_{\,W_{j}}\,U\,R^{\,\ast}\,f,\, \Gamma_{j}\,P_{\,W_{j}}\,T\,R^{\,\ast}\,f\,\right>\,\right\}\\ 
&\geq\, \sum\limits_{\,j \,\in\, J}\,v_{j}^{\,2}\,\left<\,\Lambda_{j}\,P_{\,W_{j}}\,U\,R^{\,\ast}\,f,\, \Lambda_{j}\,P_{\,W_{j}}\,T\,R^{\,\ast}\,f\,\right> \\
&\geq\; A\; \left\|\,K^{\,\ast}\,\left(\,V \,+\, W\,\right)^{\,\ast}\,f\,\right\|^{\,2}\\
&=\; A\,\left\|\,\left(\,V \,+\, W\,\right)^{\,\ast}\,K^{\,\ast}\,f\,\right\|^{\,2}\;  \;[\;\text{using}\; \,K\, \left(\,V \,+\, W \,\right) \,=\, \left(\,V \,+\, W \,\right)\, K\;]\\
&\geq\; A\,\left\|\,\left(\,V \,+\, W\,\right)^{\,-\, 1}\,\right\|^{\,-\, 2}\; \left\|\,K^{\,\ast}\,f\,\right\|^{\,2}\;  \;[\;\text{using (\ref{eqn4})}\;].
\end{align*}
Thus, \,$\Delta_{T\,U}$\, is a \,$(\,T,\,U\,)$-controlled \,$K$-$g$-fusion frame for \,$H$\, and hence by the Theorem (\ref{thm1}), it is a \,$(\,T,\,U\,)$-controlled \,$g$-atomic subspace of \,$H$\, with respect to \,$K$.\,This completes the proof.
\end{proof}

\subsection{Example}
Let \,$H \,=\, L^{\,2}\,[\,0,\,2\,\pi\,]$.\,Then \,$H$\, is a Hilbert space with respect to the inner product 
\[\left<\,x,\,y\,\right> \,=\, \int\limits_{0}^{\,2\,\pi}\,x\,(\,t\,)\,\overline{y\,(\,t\,)}\,d\,t,\; \;x,\,y \,\in\, L^{\,2}\,[\,0,\,2\,\pi\,]\]
equipped with the induced norm given by
\[\|\,x\,\|^{\,2} \,=\, \int\limits_{0}^{\,2\,\pi}\,|\,x\,(\,t\,)\,|^{\,2}\,d\,t,\; \;x \,\in\, L^{\,2}\,[\,0,\,2\,\pi\,].\]
Consider the orthonormal basis \,$\{\,u_{\,n}\,\}$, where 
\[u_{\,n}\,(\,t\,) \,=\, \frac{e^{\,i\,n\,t}}{\sqrt{2\,\pi}},\, t \,\in\, [\,0,\,2\,\pi\,],\, \,n \,\in\, \mathcal{Z}.\]
Then 
\[x \,=\, \sum\limits_{n \,=\, -\,\infty}^{\,\infty}\,x^{\prime}\,(\,n\,)\,u_{\,n},\; \;\|\,x\,\|^{\,2} \,=\, \sum\limits_{n \,=\, -\,\infty}^{\,\infty}\,|\,x^{\prime}\,(\,n\,)\,|^{\,2},\]
where 
\[x^{\prime}\,(\,n\,) \,=\, \frac{1}{\sqrt{2\,\pi}}\,\int\limits_{0}^{\,2\,\pi}\,x\,(\,t\,)\,e^{\,-\, i\,n\,t}\,d\,t.\]
Let \,$W_{n} \,=\, \overline{\,span}\,\{\,u_{\,n}\,\}$\, and $\,v_{\,n} \,=\, 1\, \,\forall\, n$.\,Now, we define the bounded linear operators \,$K \,:\, H \,\to\, H$\, by \,$K\,u_{\,1} \,=\, u_{\,1},\, \,K\,u_{\,2} \,=\, u_{\,2}$, \,$K\,u_{n} \,=\, 0$\, for \,$n \,\neq\, 1,\, 2$\, and 
\[\Lambda_{1}\,(\,x\,) \,=\, \sum\limits_{k \,=\, 1}^{m}\,x^{\prime}\,(\,k\,)\,u_{\,k},\, \;\Lambda_{j}\,(\,x\,) \,=\, 0,\, \,j \,\neq\, 1.\]
It is easy to verify that \,$K^{\,\ast}\,u_{\,1} \,=\, u_{\,1}$, \,$K^{\,\ast}\,u_{\,2} \,=\, u_{\,2}$\, and \,$K^{\,\ast}\,u_{\,n} \,=\, 0$\, for \,$n \,\neq\, 1,\, 2$.\,Consider two operators \,$T$\, and \,$U$\, on \,$H$\, defined by \,$T\,x \,=\, \alpha\,x$\, and \,$U\,x\, =\, \beta\,x$\, for \,$x \,\in\, H$\, and \,$\alpha,\, \,\beta$\, are positive real numbers.\,It is easy to see that \,$T$\, and \,$U$\, are positives.\,Now, for each \,$x \,\in\, H$, we have
\[\left\|\,K^{\,\ast}\,x\,\right\|^{\,2} \,=\, \left\|\,\sum\limits_{n \,=\, -\,\infty}^{\,\infty}\,x^{\prime}\,(\,n\,)\,K^{\,\ast}\,u_{\,n}\,\right\|^{\,2} \,=\, \left|\,x^{\prime}\,(\,1\,)\,\right|^{\,2} \,+\, \left|\,x^{\prime}\,(\,2\,)\,\right|^{\,2},\] and
\[\sum\limits_{n \,=\, -\,\infty}^{\,\infty}\,v^{\,2}_{n}\,\left<\,\Lambda_{n}\,P_{\,W_{n}}\,U\,x,\,  \Lambda_{j}\,P_{\,W_{n}}\,T\,x\,\right> \,=\, \alpha\,\beta\,\sum\limits_{n \,=\, -\,\infty}^{\,\infty}\,\left<\,\Lambda_{n}\,P_{\,W_{n}}\,x,\,  \Lambda_{n}\,P_{\,W_{n}}\,x\,\right>\]
\[=\,\alpha\,\beta\,\sum\limits_{n \,=\, -\,\infty}^{\,\infty}\,\left\|\,\Lambda_{n}\,P_{\,W_{n}}\,x\,\right\|^{\,2} \,=\, \alpha\,\beta\,\left\|\,\Lambda_{1}\,P_{\,W_{1}}\,x\,\right\|^{\,2} \,=\, \alpha\,\beta\,\left\|\,\sum\limits_{k \,=\, 1}^{m}\,x^{\prime}\,(\,k\,)\,u_{\,k}\,\right\|^{\,2}\]    
\[=\, \alpha\,\beta\,\sum\limits_{k \,=\, 1}^{m}\,|\,x^{\prime}\,(\,k\,)\,|^{\,2} \,\geq\, \alpha\,\beta\,\left\|\,K^{\,\ast}\,x\,\right\|^{\,2}.\]
Thus, for each \,$x \,\in\, H$, we have 
\[\alpha\,\beta\,\left\|\,K^{\,\ast}\,x\,\right\|^{\,2} \,\leq\, \sum\limits_{n \,=\, -\,\infty}^{\,\infty}\,v^{\,2}_{n}\,\left<\,\Lambda_{n}\,P_{\,W_{n}}\,U\,x,\,  \Lambda_{j}\,P_{\,W_{n}}\,T\,x\,\right> \,\leq\, \|\,x\,\|^{\,2}.\]
Hence, \,$\left\{\,\left(\,W_{n},\, \Lambda_{n},\, v_{n}\,\right)\,\right\}_{n \,=\, \,-\, \infty}^{\,\infty}$\, is a \,$(\,T,\,U\,)$-controlled \,$K$-$g$-fusion frame for \,$H \,=\, L^{\,2}\,[\,0,\,2\,\pi\,]$.\\  

Let \,$H$\, and \,$X$\, be two Hilbert spaces with inner products \,$\left<\,\cdot,\, \cdot\,\right>_{H}$\, and \,$\left<\,\cdot,\, \cdot\,\right>_{X}$. Then the space define by 
\[H\,\oplus\,X \,=\, \left\{\, f \,\oplus\, g \,=\, (\,f,\, g\,) \;:\; f \,\in\, H,\; g \,\in\, X\,\right\}\]
is a linear space with respect to the addition and scalar multiplication defined by
\[\left(\,f_{\,1},\, g_{\,1}\,\right) \,+\, \left(\,f_{\,2},\, g_{\,2}\,\right) \,=\, \left(\,f_{\,1} \,+\, f_{\,2},\, g_{\,1} \,+\, g_{\,2}\,\right),\;\text{and}\]
\[\lambda\,(\,f,\, g\,) \,=\, (\,\lambda\,f,\, \lambda\,g\,)\; \;\forall\; \,f,\, f_{\,1},\, f_{\,2} \,\in\, H,\;  \;g,\, g_{\,1},\, g_{\,2} \,\in\, X\;\text{and}\; \;\lambda \,\in\, \mathbb{K}.\]
Now, \,$H \,\oplus\, X$\, is an inner product space with respect to the inner product given by  
\[\left<\,(\,f \,\oplus\, g\,),\, (\,f^{\,\prime} \,\oplus\, g^{\,\prime}\,)\,\right> \,=\, \left<\,f,\, f^{\,\prime}\,\right>_{H} \,+\, \left<\,g,\, g^{\,\prime}\,\right>_{X}\; \;\forall\, f,\, f^{\,\prime} \,\in\, H\; \;\text{and}\; \;\forall\, g,\, g^{\,\prime} \,\in\, X.\]
The norm on \,$H \,\oplus\, X$\, is defined by 
\[\|\,f \,\oplus\, g\,\| \,=\, \|\,f\,\|_{H} \,+\, \|\,g\,\|_{X}\;\;  \;\forall\; \,f \,\in\, H,\;g \,\in\, X,\] where \,$\|\,\cdot\,\|_{H}$\, and \,$\|\,\cdot\,\|_{X}$\, are norms generated by \,$\left<\,\cdot,\, \cdot\,\right>_{H}$\, and \,$\left<\,\cdot,\, \cdot\,\right>_{X}$, respectively.\\The space \,$H \,\oplus\, X$\, is complete with respect to the above inner product.\,Therefore the space \,$H \,\oplus\, X$\, is a Hilbert space.\,Also, if \,$R \,\in\, \mathcal{B}\,(\,H,\, Z\,),\, V \,\in\, \mathcal{B}\,(\,X,\, Y)$, where \,$Z,\, Y$\, are Hilbert spaces then for all \,$f \,\in\, H$\, and \,$g \,\in\, X$\, we define
\[R\,\oplus\,V \,\in\, \mathcal{B}\,\left(\,H\,\oplus\,X,\, Z\,\oplus\,Y \right)\; \;\text{by}\; \left(\,R\,\oplus\,V \,\right)\,(\,f \,\oplus\, g\,) \,=\, (\,R\,f,\, V\,g\,)\]
Now, we state a theorem, for proof of this theorem, one can go through the paper \cite{PG}.
\begin{theorem}\cite{PG}
Suppose  \,$R,\, R^{\,\prime} \,\in\, \mathcal{B}\,(\,H,\, Z\,)$\, and \,$V,\, V^{\,\prime} \,\in\, \mathcal{B}\,(\,X,\, Y)$.\;Then 
\begin{description}
\item[$(i)$]$\left(\,R \,\oplus\, V\,\right)\,\left(\,R^{\,\prime} \,\oplus\, V^{\,\prime}\,\right) \,=\, \left(\,R\,R^{\,\prime} \,\oplus\, V\,V^{\,\prime}\,\right)$.
\item[$(ii)$]\, $\left(\,R \,\oplus\,V \,\right)^{\,\ast} \,=\, R^{\,\ast} \,\oplus\, V^{\,\ast}$.
\item[$(iii)$] If \,$R\; \;\text{and}\; \,V$\, are invertible, then  \,$\left(\,R\,\oplus\,V\,\right)$\, is invertible and moreover we have \,$\left(\,R\,\oplus\,V\,\right)^{\,-\, 1} \,=\, R^{\,-\, 1} \,\oplus\, V^{\,-\, 1}$.
\item[$(iv)$]\,$I_{H} \,\oplus\, I_{K} \,=\, I_{H \,\oplus\, K}$, where \,$I_{H}$, \,$I_{K}$\, and \,$I_{H \,\oplus\, K}$\, are identity operators on \,$H$, \,$K$\, and \,$H \,\oplus\, K$, respectively.
\item[$(v)$]\,$P_{\,M \,\oplus\, N} \,=\, P_{M} \,\oplus\, P_{N}$, where \,$P_{\,M},\, P_{\,N}\; \;\text{and}\; \;P_{\,M\,\oplus\,\,N}$\; are the orthogonal projection onto the closed subspaces  \;$M \,\subset\, H,\, N \,\subset\, X \;\text{and}\; \;M\,\oplus\,\,N \,\subset\, H\,\oplus\,X$, respectively.
\end{description}
\end{theorem}

We now construct controlled \,$K$-$g$-fusion frame for the Hilbert space \,$H \,\oplus\, X$, using the controlled \,$K$-$g$-fusion frames of the Hilbert spaces \,$H$\, and \,$X$.\,Assuming for each \,$j \,\in\, J,\, \,W_{j}\,\oplus\,V_{j}$\, are the closed subspaces of \,$H\,\oplus\,X$, \,$\Gamma_{j} \,\in\, \mathcal{B}\,(\,X,\, X_{j})$, where \,$\{\,X_{j}\,\}_{j \,\in\, J}$\, are the collection of Hilbert spaces and \,$\Lambda_{j}\,\oplus\,\Gamma_{j} \,\in\, \mathcal{B}\,\left(\,H\,\oplus\,X,\, H_{j} \,\oplus\, X_{j} \right)$.

\begin{theorem}
Let \,$\Lambda_{T\,T_{\,1}} \,=\, \left\{\,\left(\,W_{j},\, \Lambda_{j},\, v_{j}\,\right)\,\right\}_{j \,\in\, J}$\, be a \,$(\,T,\,T_{\,1}\,)$-controlled \,$K_{1}$-$g$-fusion frame for \,$H$\, with bounds \,$A,\, B$\, and \,$\Gamma_{U\,U_{\,1}} \,=\, \left\{\,\left(\,V_{j},\, \Gamma_{j},\, v_{j}\,\right)\,\right\}_{j \,\in\, J}$\, be a \,$(\,U,\,U_{\,1}\,)$-controlled \,$K_{2}$-$g$-fusion frame for \,$X$\, with bounds \,$C,\, D$, where \,$K_{1},\, K_{2} \,\in\, \mathcal{B}\,(\,H\,)$\, and \,$T_{1},\, U_{1} \,\in\, \mathcal{G}\,\mathcal{B}\,(\,H\,)$.\,Then \,$\Lambda_{T\,T_{\,1}} \,\oplus\, \Gamma_{U\,U_{\,1}} \,=\, \left\{\,\left(\,W_{j}\,\oplus\,V_{j},\; \Lambda_{j}\,\oplus\,\Gamma_{j},\; v_{j}\,\right)\,\right\}_{j \,\in\, J}$\, is a \,$\left(\,T \,\oplus\, U,\, T_{\,1} \,\oplus\, U_{\,1}\,\right)$-controlled $K_{\,1} \,\oplus\, K_{\,2}$-$g$-fusion frame for \,$H\,\oplus\,X$\, with bounds \,$\min\,\{\,A,\, C\,\},\; \max\,\{\,B,\, D\,\}$.\,Furthermore, if \,$S_{C},\, S_{C^{\prime}}$\, and \,$S$\, are frame operators for \,$\Lambda_{T\,T_{\,1}},\, \Gamma_{U\,U_{\,1}}$\, and \,$\Lambda_{T\,T_{\,1}} \,\oplus\, \Gamma_{U\,U_{\,1}}$, respectively then we have \,$S \,=\, S_{C} \,\oplus\, S_{C^{\prime}}$. 
\end{theorem}

\begin{proof}
For each \,$f \,\oplus\, g \,\in\, H\,\oplus\,X$, we have
\begin{align*}
&\sum\limits_{\,j \,\in\, J}\,v_{j}^{\,2}\,\left<\,\left(\Lambda_{j}\,\oplus\,\Gamma_{j}\right)P_{\,W_{j}\,\oplus\,V_{j}}\left(T_{1} \,\oplus\, U_{1}\right)(f \,\oplus\, g),\, \left(\Lambda_{j}\,\oplus\,\Gamma_{j}\right)P_{\,W_{j}\,\oplus\,V_{j}}\left(T \,\oplus\, U\right)(f \,\oplus\, g)\,\right>\\
&=\,\sum\limits_{\,j \,\in\, J}\,v_{j}^{\,2}\,\left<\,\left(\,\Lambda_{j}\,\oplus\,\Gamma_{j}\,\right)\, P_{\,W_{j}\,\oplus\,V_{j}}\,\left(\,T_{\,1}\,f \,\oplus\, U_{\,1}\,g\,\right),\, \left(\,\Lambda_{j}\,\oplus\,\Gamma_{j}\,\right)\, P_{\,W_{j}\,\oplus\,V_{j}}\,\left(\,T\,f \,\oplus\, U\,g\,\right)\,\right>\\
&=\,\sum\limits_{\,j \,\in\, J}\,v_{j}^{\,2}\,\left<\,\left(\,\Lambda_{j}\,\oplus\,\Gamma_{j}\,\right)\,\left(\,P_{\,W_{j}}\,T_{\,1}\,f \,\oplus\, P_{\,V_{j}}\,U_{\,1}\,g\,\right),\, \left(\,\Lambda_{j}\,\oplus\,\Gamma_{j}\,\right)\,\left(\,P_{\,W_{j}}\,T\,f \,\oplus\, P_{\,V_{j}}\,U\,g\,\right)\,\right>\\
&=\,\sum\limits_{\,j \,\in\, J}\,v_{j}^{\,2}\,\left<\,\left(\,\Lambda_{j}\,P_{\,W_{j}}\,T_{\,1}\,f \,\oplus\, \Gamma_{j}\,P_{\,V_{j}}\,U_{\,1}\,g\,\right),\, \left(\,\Lambda_{j}P_{\,W_{j}}\,T\,f \,\oplus\, \Gamma_{j}\,P_{\,V_{j}}\,U\,g\,\right)\,\right>\\
&=\, \sum\limits_{\,j \,\in\, J}\,v_{j}^{\,2}\,\left<\,\Lambda_{j}\,P_{\,W_{j}}\,T_{\,1}\,f,\, \Lambda_{j}P_{\,W_{j}}\,T\,f \,\right>_{H} \,+\, \sum\limits_{\,j \,\in\, J}\,v_{j}^{\,2}\,\left<\,\Gamma_{j}\,P_{\,V_{j}}\,U_{\,1}\,g,\, \Gamma_{j}\,P_{\,V_{j}}\,U\,g\,\right>_{X}\\ 
&\leq\, B\,\|\,f\,\|_{H}^{\,2} \,+\, D\,\|\,g\,\|_{X}^{\,2} \,\leq\, \max\,\{\,B,\, D\,\}\, \left(\,\|\,f\,\|_{H}^{\,2} \,+\, \|\,g\,\|_{X}^{\,2}\,\right)\\
& \,=\, \max\,\{\,B,\, D\,\}\,\|\,f \,\oplus\, g\,\|^{\,2}.
\end{align*}
Similarly, for each \,$f \,\oplus\, g \,\in\, H\,\oplus\,X$, it can be shown that
\begin{align*}
&\sum\limits_{\,j \,\in\, J}\,v_{j}^{\,2}\left<\,\left(\Lambda_{j}\,\oplus\,\Gamma_{j}\right)P_{\,W_{j}\,\oplus\,V_{j}}\left(T_{1} \,\oplus\, U_{1}\right)(f \,\oplus\, g),\, \left(\Lambda_{j}\,\oplus\,\Gamma_{j}\right)P_{\,W_{j}\,\oplus\,V_{j}}\left(T \,\oplus\, U\right)(f \,\oplus\, g)\,\right>\\
&\hspace{1cm}\geq\, \min\,\{\,A,\, C\,\}\,\left(\,\left\|\,K_{1}^{\,\ast}\,f\,\right\|_{H}^{\,2} \,+\, \left\|\,K_{2}^{\,\ast}\,g\,\right\|_{X}^{\,2}\,\right) \,=\, \min\,\{\,A,\, C\,\}\,\left\|\,K^{\,\ast}_{\,1}\,f \,\oplus\, K^{\,\ast}_{\,2}\,g\,\right\|^{\,2}\\
&\hspace{1cm} \,=\, \min\,\{\,A,\, C\,\}\,\left\|\,(\,K_{\,1} \,\oplus\, K_{\,2}\,)^{\,\ast}\,(\,f \,\oplus\, g\,)\,\right\|^{\,2}.
\end{align*}
Thus, \,$\Lambda_{T\,T_{\,1}} \,\oplus\, \Gamma_{U\,U_{\,1}}$\, is a \,$\left(\,T \,\oplus\, U,\, T_{\,1} \,\oplus\, U_{\,1}\,\right)$-controlled $K_{\,1} \,\oplus\, K_{\,2}$-$g$-fusion frame for \,$H\,\oplus\,X$\, with bounds \,$\min\,\{\,A,\, C\,\}$\, and \,$\max\,\{\,B,\, D\,\}$.\\

Furthermore, for \,$f \,\oplus\, g \,\in\, H\,\oplus\,X$, we have 
\begin{align*}
&S\,(\,f \,\oplus\, g\,)\\
& \,=\, \sum\limits_{\,j \,\in\, J}\,v_{j}^{\,2}\left(\,T \,\oplus\, U\,\right)^{\,\ast}\,P_{\,W_{j}\,\oplus\,V_{j}}\left(\,\Lambda_{j}\,\oplus\,\Gamma_{j}\,\right)^{\,\ast}\left(\,\Lambda_{j}\,\oplus\,\Gamma_{j}\,\right)P_{\,W_{j}\,\oplus\,V_{j}}\left(\,T_{\,1} \,\oplus\, U_{\,1}\,\right)(f \,\oplus\, g)\\ 
& =\; \sum\limits_{\,j \,\in\, J}\,v_{j}^{\,2}\,\left(\,T \,\oplus\, U\,\right)^{\,\ast}\,P_{\,W_{j}\,\oplus\,V_{j}}\, \left(\,\Lambda_{j}\,\oplus\,\Gamma_{j}\,\right)^{\,\ast}\, \left(\,\Lambda_{j}\,P_{\,W_{j}}\,T_{\,1}\,f \,\oplus\, \Gamma_{j}\,P_{\,V_{j}}\,U_{\,1}\,g\,\right)\\
& =\, \sum\limits_{\,j \,\in\, J}\,v_{j}^{\,2}\,\left(\,T \,\oplus\, U\,\right)^{\,\ast}\,P_{\,W_{j}\,\oplus\,V_{j}}\, \left(\,\Lambda_{j}^{\,\ast}\, \Lambda_{j}\, P_{\,W_{j}}\,T_{\,1}\,f \,\oplus\, \Gamma_{j}^{\,\ast}\, \Gamma_{j}\, P_{\,V_{j}}\,U_{\,1}\,g\,\right)\\
&=\; \sum\limits_{\,j \,\in\, J}\,v_{j}^{\,2}\; \left(\,T^{\,\ast}\,P_{\,W_{j}}\,\Lambda_{j}^{\,\ast}\, \Lambda_{j}\, P_{\,W_{j}}\,T_{\,1}\,f \,\oplus\, U^{\,\ast}\,P_{\,V_{j}}\, \Gamma_{j}^{\,\ast}\, \Gamma_{j}\,P_{\,V_{j}}\,U_{\,1}\,g\,\right)\hspace{.1cm}\\
&=\; \left(\,\sum\limits_{\,j \,\in\, J}\,v_{j}^{\,2}\,T^{\,\ast}\,P_{\,W_{j}}\, \Lambda_{j}^{\,\ast}\, \Lambda_{j}\, P_{\,W_{j}}\,T_{\,1}\,f\,\right) \,\oplus\, \left(\,\sum\limits_{\,j \,\in\, J}\,v_{j}^{\,2}\,U^{\,\ast}\,P_{\,V_{j}}\,\Gamma_{j}^{\,\ast}\, \Gamma_{j}\, P_{\,V_{j}}\,U_{\,1}\,g\,\right)\\
&=\; \left(\,S_{C}\,f \,\oplus\, S_{C^{\,\prime}}\,g\,\right) \,=\, \left(\,S_{C} \,\oplus\, S_{C^{\,\prime}}\,\right)\,(\,f \,\oplus\, g\,).
\end{align*} 
Hence, \,$S \,=\, S_{C} \,\oplus\, S_{C^{\,\prime}}$.\,This completes the proof.    
\end{proof}

In the following theorem, we give another way to construct controlled $K_{\,1} \,\oplus\, K_{\,2}$-$g$-fusion frame for Hilbert space \,$H\,\oplus\,X$. 

\begin{theorem}
Let \,$\Lambda_{T\,T_{\,1}} \,=\, \left\{\,\left(\,W_{j},\, \Lambda_{j},\, v_{j}\,\right)\,\right\}_{j \,\in\, J}$\, be a \,$(\,T,\,T_{\,1}\,)$-controlled \,$K_{1}$-$g$-fusion frame for \,$H$\, with bounds \,$A,\, B$\, and \,$\Gamma_{U\,U_{\,1}} \,=\, \left\{\,\left(\,V_{j},\, \Gamma_{j},\, v_{j}\,\right)\,\right\}_{j \,\in\, J}$\, be a \,$(\,U,\,U_{\,1}\,)$-controlled \,$K_{2}$-$g$-fusion frame for \,$X$\, with bounds \,$C,\, D$\, having their corresponding frame operators \,$S_{C}$\, and \,$S_{C^{\,\prime}}$, where \,$K_{1},\, K_{2} \,\in\, \mathcal{B}\,(\,H\,)$\, and \,$T_{1},\, U_{1} \,\in\, \mathcal{G}\,\mathcal{B}\,(\,H\,)$.\,Suppose 
\begin{description}
\item[$(i)$]\,$W \,\in\, \mathcal{B}\,(\,H\,)$\, and \,$V \,\in\, \mathcal{B}\,(\,X\,)$\, are invertible operators on \,$H$\, and \,$X$, respectively.
\item[$(ii)$]\,$W^{\,\ast}$\, commutes with \,$T,\, T_{1}$\, and \,$V^{\,\ast}$\, commutes with \,$U,\, U_{1}$, respectively.
\item[$(iii)$] \,$K_{1}$\, and \,$K_{2}$\, commutes with \,$W$\, and \,$V$, respectively.
\end{description}
Then  
\[\Delta \,=\, \left\{\,\left(\,\left(\,W \,\oplus\, V\,\right)\,\left(\,W_{j}\,\oplus\,V_{j}\,\right),\; \left(\,\Lambda_{j}\,\oplus\,\Gamma_{j}\,\right)\,P_{\,W_{j}\,\oplus\,V_{j}}\,\left(\,W^{\,\ast} \,\oplus\, V^{\,\ast}\,\right),\; v_{j}\,\right)\,\right\}_{j \,\in\, J}\]
 is a \,$\left(\,T \,\oplus\, U,\, T_{\,1} \,\oplus\, U_{\,1}\,\right)$-controlled $K_{\,1} \,\oplus\, K_{\,2}$-$g$-fusion frame for \,$H\,\oplus\,X$\, with the corresponding frame operator \,$(\,W \,\oplus\, V\,)\,\left(\,S_{C} \,\oplus\, S_{C^{\,\prime}}\,\right)\,\left(\,W \,\oplus\, V\,\right)^{\,\ast}$.   
\end{theorem}

\begin{proof}
For \,$j \,\in\, J$, we take 
\[\Delta_{j} \,=\, \left(\,\Lambda_{j}\,\oplus\,\Gamma_{j}\,\right)\,P_{\,W_{j}\,\oplus\,V_{j}}\,\left(\,W^{\,\ast} \,\oplus\, V^{\,\ast}\,\right),\, \,X_{j} \,=\, \left(\,W \,\oplus\, V\,\right)\,\left(\,W_{j}\,\oplus\,V_{j}\,\right).\]
Then by Theorem \ref{th1.01}, we have 
\begin{align*}
&\Delta_{j}\,P_{X_{j}}\,(\,f \,\oplus\, g\,) \\
&\,=\, \left(\,\Lambda_{j}\,\oplus\,\Gamma_{j}\,\right)\,\left(\,P_{\,W_{j}} \,\oplus\, P_{\,V_{j}}\,\right)\,\left(\,W^{\,\ast} \,\oplus\, V^{\,\ast}\,\right)\,\left(\,P_{\,W\,W_{j}} \,\oplus\, P_{\,V\,V_{j}}\,\right)\,(\,f \,\oplus\, g\,)\\
&=\,\Lambda_{j}\,P_{\,W_{j}}\,W^{\,\ast}\,P_{\,W\,W_{j}}\,f \,\oplus\, \Gamma_{j}\,P_{\,V_{j}}\,V^{\,\ast}\,P_{\,V\,V_{j}}\,g\\
&\,=\,\Lambda_{j}\,P_{\,W_{j}}\,W^{\,\ast}\,f \,\oplus\, \Gamma_{j}\,P_{\,V_{j}}\,V^{\,\ast}\,g.    
\end{align*} 
and
\begin{align*}
&P_{X_{j}}\,\Delta_{j}^{\,\ast}\,(\,f \,\oplus\, g\,)\\
&=\,\left(\,P_{\,W\,W_{j}} \,\oplus\, P_{\,V\,V_{j}}\,\right)\,(\,W \,\oplus\, V\,)\,\left(\,P_{\,W_{j}} \,\oplus\, P_{\,V_{j}}\,\right)\,\left(\,\Lambda^{\,\ast}_{j} \,\oplus\, \Gamma^{\,\ast}_{j}\,\right)\,(\,f \,\oplus\, g\,)\\
&=\,P_{\,W\,W_{j}}\,W\,P_{\,W_{j}}\,\Lambda^{\,\ast}_{j}\,f \,\oplus\, P_{\,V\,V_{j}}\,V\,P_{\,V_{j}}\,\Gamma^{\,\ast}_{j}\,g\\
&=\, \left(\,P_{\,W_{j}}\,W^{\,\ast}\,P_{\,W\,W_{j}}\,\right)^{\,\ast}\,\Lambda^{\,\ast}_{j}\,f \,\oplus\, \left(\,P_{\,V_{j}}\,V^{\,\ast}\,P_{\,V\,V_{j}}\,\right)^{\,\ast}\,\Gamma^{\,\ast}_{j}\,g\\
&=\, W\,P_{\,W_{j}}\,\Lambda^{\,\ast}_{j}\,f \,\oplus\, V\,P_{\,V_{j}}\,\Gamma^{\,\ast}_{j}\,g.  
\end{align*}
Therefore,
\begin{align}
&P_{X_{j}}\,\Delta_{j}^{\,\ast}\,\Delta_{j}\,P_{X_{j}}\,(\,f \,\oplus\, g\,)\nonumber\\
&=\,W\,P_{\,W_{j}}\,\Lambda^{\,\ast}_{j}\,\Lambda_{j}\,P_{\,W_{j}}\,W^{\,\ast}\,f \,\oplus\, V\,P_{\,V_{j}}\,\Gamma^{\,\ast}_{j}\,\Gamma_{j}\,P_{\,V_{j}}\,V^{\,\ast}\,g.\label{eqnn1} 
\end{align}
Now, for each \,$f \,\oplus\, g \,\in\, H\,\oplus\,X$, we have
\begin{align}
&\sum\limits_{\,j \,\in\, J}\,v_{j}^{\,2}\,\left<\,\Delta_{j}\,P_{X_{j}}\,\left(T_{1} \,\oplus\, U_{1}\right)\,(\,f \,\oplus\, g\,), \Delta_{j}\,P_{X_{j}}\,\left(\,T \,\oplus\, U\,\right)\,(\,f \,\oplus\, g\,)\,\right>\nonumber\\
&=\, \sum\limits_{\,j \,\in\, J}\,v_{j}^{\,2}\,\left<\,\Delta_{j}\,P_{X_{j}}\,\left(\,T_{1}\,f \,\oplus\, U_{1}\,g\right), \Delta_{j}\,P_{X_{j}}\,\left(\,T\,f \,\oplus\, U\,g\,\right)\,\right>\nonumber\\
&=\, \sum\limits_{\,j \,\in\, J}\,v_{j}^{\,2}\,\left<\,\Lambda_{j}\,P_{\,W_{j}}\,W^{\,\ast}\,T_{1}\,f \,\oplus\, \Gamma_{j}\,P_{\,V_{j}}\,V^{\,\ast}\,U_{1}\,g,\, \Lambda_{j}\,P_{\,W_{j}}\,W^{\,\ast}\,T\,f \,\oplus\, \Gamma_{j}\,P_{\,V_{j}}\,V^{\,\ast}\,U\,g\right>\nonumber\\
&=\, \sum\limits_{\,j \,\in\, J}\,v_{j}^{\,2}\,\left<\,\Lambda_{j}\,P_{\,W_{j}}\,W^{\,\ast}\,T_{\,1}\,f,\, \Lambda_{j}P_{\,W_{j}}\,W^{\,\ast}\,T\,f \,\right>_{H} \nonumber\\
&\hspace{1cm}\,+\, \sum\limits_{\,j \,\in\, J}\,v_{j}^{\,2}\,\left<\,\Gamma_{j}\,P_{\,V_{j}}\,V^{\,\ast}\,U_{\,1}\,g,\, \Gamma_{j}\,P_{\,V_{j}}\,V^{\,\ast}\,U\,g\,\right>_{X}\nonumber\\
&=\, \sum\limits_{\,j \,\in\, J}\,v_{j}^{\,2}\,\left<\,\Lambda_{j}\,P_{\,W_{j}}\,T_{\,1}\,W^{\,\ast}\,f,\, \Lambda_{j}P_{\,W_{j}}\,T\,W^{\,\ast}\,f \,\right>_{H} \nonumber\\
&\hspace{1cm}\,+\, \sum\limits_{\,j \,\in\, J}\,v_{j}^{\,2}\,\left<\,\Gamma_{j}\,P_{\,V_{j}}\,U_{\,1}\,V^{\,\ast}\,g,\, \Gamma_{j}\,P_{\,V_{j}}\,U\,V^{\,\ast}\,g\,\right>_{X}\label{eq1.1}\\
&\geq\,A_{1}\,\left\|\,K_{1}^{\,\ast}\,W^{\,\ast}\,f\,\right\|_{H}^{\,2} \,+\, A_{2}\,\left\|\,K_{2}^{\,\ast}\,V^{\,\ast}\,g\,\right\|_{X}^{\,2} \,=\, A_{1}\,\left\|\,W^{\,\ast}\,K_{1}^{\,\ast}\,f\,\right\|_{H}^{\,2} \,+\, A_{2}\,\left\|\,V^{\,\ast}\,K_{2}^{\,\ast}\,g\,\right\|_{X}^{\,2}\nonumber\\
&\geq\, \dfrac{A_{1}}{\left\|\,W^{\,-\, 1}\,\right\|^{\,2}}\,\left\|\,K_{1}^{\,\ast}\,f\,\right\|_{H}^{\,2} \,+\, \dfrac{A_{2}}{\left\|\,V^{\,-\, 1}\,\right\|^{\,2}}\,\left\|\,K_{2}^{\,\ast}\,g\,\right\|_{X}^{\,2}\; \;[\;\text{since $W,\,V$ are invertible}\;]\nonumber\\
&\geq\, C\,\left\|\,K^{\,\ast}_{\,1}\,f \,\oplus\, K^{\,\ast}_{\,2}\,g\,\right\|^{\,2} \,=\, C\,\left\|\,(\,K_{\,1} \,\oplus\, K_{\,2}\,)^{\,\ast}\,(\,f \,\oplus\, g\,)\,\right\|^{\,2},\nonumber    
\end{align}
where \,$C \,=\, \min\,\left(\,\dfrac{A_{1}}{\left\|\,W^{\,-\, 1}\,\right\|^{\,2}},\, \dfrac{A_{2}}{\left\|\,V^{\,-\, 1}\,\right\|^{\,2}}\,\right)$.\\On the other hand, for each \,$f \,\oplus\, g \,\in\, H\,\oplus\,X$, from (\ref{eq1.1}), we have
\begin{align*}
&\sum\limits_{\,j \,\in\, J}\,v_{j}^{\,2}\,\left<\,\Delta_{j}\,P_{X_{j}}\,\left(T_{1} \,\oplus\, U_{1}\right)\,(\,f \,\oplus\, g\,), \Delta_{j}\,P_{X_{j}}\,\left(\,T \,\oplus\, U\,\right)\,(\,f \,\oplus\, g\,)\,\right>\\
&\leq\, B_{1}\,\left\|\,W^{\,\ast}\,f\,\right\|_{H}^{\,2} \,+\, B_{2}\,\left\|\,V^{\,\ast}\,g\,\right\|_{X}^{\,2} \,\leq\, \max\,\left(\,B_{1}\,\|\,W\,\|^{\,2},\, B_{2}\,\|\,V\,\|^{\,2}\,\right)\,\|\,f \,\oplus\, g\,\|^{\,2}.
\end{align*}
Thus, \,$\Delta$\, is a \,$\left(\,T \,\oplus\, U,\, T_{\,1} \,\oplus\, U_{\,1}\,\right)$-controlled $K_{\,1} \,\oplus\, K_{\,2}$-$g$-fusion frame for \,$H\,\oplus\,X$.\\

Furthermore, for each \,$f \,\oplus\, g \,\in\, H\,\oplus\,X$, using (\ref{eqnn1}), we have
\begin{align*}
&\sum\limits_{\,j \,\in\, J}\,v_{j}^{\,2}\,\left(\,T \,\oplus\, U\,\right)^{\,\ast}\,P_{X_{j}}\,\Delta_{j}^{\,\ast}\,\Delta_{j}\,P_{X_{j}}\,\left(\,T_{\,1} \,\oplus\, U_{\,1}\,\right)\,(\,f \,\oplus\, g\,)\\
&=\, \sum\limits_{\,j \,\in\, J}\,v_{j}^{\,2}\,\left(\,T^{\,\ast} \,\oplus\, U^{\,\ast}\,\right)\,P_{X_{j}}\,\Delta_{j}^{\,\ast}\,\Delta_{j}\,P_{X_{j}}\,\left(\,T_{\,1}\,f \,\oplus\, U_{\,1}\,g\,\right)\\
&=\, \sum\limits_{\,j \,\in\, J}\,v_{j}^{\,2}\,\left(\,T^{\,\ast}\,W\,P_{\,W_{j}}\,\Lambda^{\,\ast}_{j}\,\Lambda_{j}\,P_{\,W_{j}}\,W^{\,\ast}\,T_{\,1}\,f\,\right) \,\oplus\, \left(\,U^{\,\ast}\,V\,P_{\,V_{j}}\,\Gamma^{\,\ast}_{j}\,\Gamma_{j}\,P_{\,V_{j}}\,V^{\,\ast}\,U_{\,1}\,g\,\right)\\
&=\,\sum\limits_{\,j \,\in\, J}\,v_{j}^{\,2}\,W\,T^{\,\ast}\,P_{\,W_{j}}\,\Lambda^{\,\ast}_{j}\,\Lambda_{j}\,P_{\,W_{j}}\,T_{\,1}\,W^{\,\ast}\,f \,\oplus\, \sum\limits_{\,j \,\in\, J}\,v_{j}^{\,2}\,V\,U^{\,\ast}\,P_{\,V_{j}}\,\Gamma^{\,\ast}_{j}\,\Gamma_{j}\,P_{\,V_{j}}\,U_{\,1}\,V^{\,\ast}\,g\\
&=\, \left(\,W\,S_{C}\,W^{\,\ast} \,\oplus\, V\,S_{C^{\,\prime}}\,V^{\,\ast}\,\right)\,(\,f \,\oplus\, g\,)\\
& \,=\, (\,W \,\oplus\, V\,)\,\left(\,S_{C} \,\oplus\, S_{C^{\,\prime}}\,\right)\,\left(\,W \,\oplus\, V\,\right)^{\,\ast}\,(\,f \,\oplus\, g\,).    
\end{align*} 
Thus, the corresponding frame operator is \,$(\,W \,\oplus\, V\,)\,\left(\,S_{C} \,\oplus\, S_{C^{\,\prime}}\,\right)\,\left(\,W \,\oplus\, V\,\right)^{\,\ast}$. This completes the proof.
\end{proof}

\section{Frame operator for a pair of controlled $g$-fusion Bessel sequences}

\smallskip\hspace{.6 cm} In this section, we introduce the concept of the frame operator for a pair of controlled \,$g$-fusion Bessel sequences and establish some properties relative to frame operator.\,First, we discuss about resolutions of the identity operator on \,$H$\, using the theory of controlled $g$-fusion frame for Hilbert space.

Let \,$\Lambda_{T\,U}$\, be a \,$(\,T,\,U\,)$-controlled $g$-fusion frame for \,$H$\, with the corresponding frame operator \,$S_{C}$.\,From reconstruction formula, for \,$f \,\in\, H$, we have
\begin{align*}
f &\,=\, S_{C}\,S^{\,-\, 1}_{C}\,f \,=\, S^{\,-\, 1}_{C}\,S_{C}\,f\\
&=\, \sum\limits_{\,j \,\in\, J}\,v_{j}^{\,2}\,T^{\,\ast}\,P_{\,W_{j}}\, \Lambda_{j}^{\,\ast}\, \Lambda_{j}\, P_{\,W_{j}}\,U\,S^{\,-\, 1}_{C}\,f \,=\, \sum\limits_{\,j \,\in\, J}\,v_{j}^{\,2}\,S^{\,-\, 1}_{C}\,T^{\,\ast}\,P_{\,W_{j}}\,\Lambda_{j}^{\,\ast}\, \Lambda_{j}\, P_{\,W_{j}}\,U\,f  
\end{align*}
This shows that the families of bounded operators \,$\left\{\,v_{j}^{\,2}\,T^{\,\ast}\,P_{\,W_{j}}\, \Lambda_{j}^{\,\ast}\,\Lambda_{j}\, P_{\,W_{j}}\,U\,S^{\,-\, 1}_{C}\,\right\}_{j \,\in\, J}$\, and \,$\left\{\,v_{j}^{\,2}\,S^{\,-\, 1}_{C}\,T^{\,\ast}\,P_{\,W_{j}}\,\Lambda_{j}^{\,\ast}\,\Lambda_{j}\, P_{\,W_{j}}\,U\,\right\}_{j \,\in\, J}$\, are resolution of the identity operator on \,$H$.  

\begin{theorem}
Let \,$\Lambda_{T\,U}$\, be a \,$(\,T,\,U\,)$-controlled $g$-fusion frame for \,$H$\, with frame bounds \,$A,\,B$\, and \,$S_{C}$\, be its corresponding frame operator.\,Assume that \,$S^{\,-\, 1}_{C}$\, commutes with \,$T$\, and \,$U$.\,Then the family \,$\left\{\,v_{j}^{\,2}\,T^{\,\ast}\,P_{\,W_{j}}\,\Lambda_{j}^{\,\ast}\,T_{j}\,U\,\right\}_{j \,\in\, J}$\, is a resolution of the identity operator on \,$H$, where \,$T_{j} \,=\, \Lambda_{j}\, P_{\,W_{j}}\,S^{\,-\, 1}_{C},\, \,j \,\in\, J$.\,Furthermore, for each \,$f \,\in\, H$, we have 
\[\dfrac{A}{B^{\,2}}\,\|\,f\,\|^{\,2} \,\leq\, \sum\limits_{\,j \,\in\, J}\,v_{j}^{\,2}\,\left<\,T_{j}\,U\,f,\, T_{j}\,T\,f\,\right> \,\leq\, \dfrac{B}{A^{\,2}}\,\|\,f\,\|^{\,2}.\] 
\end{theorem}

\begin{proof}
For \,$f \,\in\, H$, we have the reconstruction formula for \,$\Lambda_{T\,U}$:
\begin{align*}
f &\,=\, \sum\limits_{\,j \,\in\, J}\,v_{j}^{\,2}\,T^{\,\ast}\,P_{\,W_{j}}\, \Lambda_{j}^{\,\ast}\, \Lambda_{j}\, P_{\,W_{j}}\,U\,S^{\,-\, 1}_{C}\,f\\
&\,=\, \sum\limits_{\,j \,\in\, J}\,v_{j}^{\,2}\,T^{\,\ast}\,P_{\,W_{j}}\, \Lambda_{j}^{\,\ast}\, \Lambda_{j}\, P_{\,W_{j}}\,S^{\,-\, 1}_{C}\,U\,f \,=\, \sum\limits_{\,j \,\in\, J}\,v_{j}^{\,2}\,T^{\,\ast}\,P_{\,W_{j}}\, \Lambda_{j}^{\,\ast}\,T_{j}\,U\,f.
\end{align*}
Thus, \,$\left\{\,v_{j}^{\,2}\,T^{\,\ast}\,P_{\,W_{j}}\,\Lambda_{j}^{\,\ast}\,T_{j}\,U\,\right\}_{j \,\in\, J}$\, is a resolution of the identity operator on \,$H$.\,Since \,$\Lambda_{T\,U}$\, is a \,$(\,T,\,U\,)$-controlled $g$-fusion frame for \,$H$\, with frame bounds \,$A$\, and \,$B$, for each \,$f \,\in\, H$, we have
\begin{align*}
\sum\limits_{\,j \,\in\, J}\,v_{j}^{\,2}\,\left<\,T_{j}\,U\,f,\, T_{j}\,T\,f\,\right>& \,=\, \sum\limits_{\,j \,\in\, J}\,v_{j}^{\,2}\,\left<\,\Lambda_{j}\, P_{\,W_{j}}\,S^{\,-\, 1}_{C}\,U\,f,\, \Lambda_{j}\, P_{\,W_{j}}\,S^{\,-\, 1}_{C}\,T\,f\,\right>\\
&=\, \sum\limits_{\,j \,\in\, J}\,v_{j}^{\,2}\,\left<\,\Lambda_{j}\,P_{\,W_{j}}\,U\,S^{\,-\, 1}_{C}\,f,\, \Lambda_{j}\,P_{\,W_{j}}\,T\,S^{\,-\, 1}_{C}\,f\,\right>\\
&\leq\, B\,\left\|\,S^{\,-\, 1}_{C}\,f\,\right\|^{\,2} \,\leq\, B\,\left\|\,S^{\,-\, 1}_{C}\,\right\|^{\,2}\,\|\,f\,\|^{\,2} \,\leq\, \dfrac{B}{A^{\,2}}\,\|\,f\,\|^{\,2}.  
\end{align*}
On the other hand, for each \,$f \,\in\, H$, we have 
\[\sum\limits_{\,j \,\in\, J}\,v_{j}^{\,2}\,\left<\,T_{j}\,U\,f,\, T_{j}\,T\,f\,\right> \,\geq\, A\,\left\|\,S^{\,-\, 1}_{C}\,f\,\right\|^{\,2} \,\geq\, \dfrac{A}{B^{\,2}}\,\|\,f\,\|^{\,2}.\]
This completes the proof.   
\end{proof}

Next we will see that a controlled $g$-fusion Bessel sequence becomes a controlled $g$-fusion frame by using a resolution of the identity operator on \,$H$.       

\begin{theorem}
Let \,$\Lambda_{T} \,=\, \left\{\,\left(\,W_{j},\, \Lambda_{j},\, v_{j}\,\right)\,\right\}_{j \,\in\, J}$\, be a \,$(\,T,\,T\,)$-controlled $g$-fusion Bessel sequence in \,$H$\, with bound \,$B$.\,Then \,$\Lambda_{T}$\, is a \,$(\,U,\,U\,)$-controlled $g$-fusion frame for  \,$H$\, if \,$\left\{\,v_{j}^{\,2}\,T^{\,\ast}\,P_{\,W_{j}}\,\Lambda_{j}^{\,\ast}\,\Lambda_{j}\,P_{\,W_{j}}\,U\,\right\}_{j \,\in\, J}$\, is a resolution of the identity operator on \,$H$.      
\end{theorem}

\begin{proof}
Since \,$\left\{\,v_{j}^{\,2}\,T^{\,\ast}\,P_{\,W_{j}}\,\Lambda_{j}^{\,\ast}\,\Lambda_{j}\,P_{\,W_{j}}\,U\,\right\}_{j \,\in\, J}$\, is a resolution of the identity operator on \,$H$, for \,$f \,\in\, H$, we have
\[f \,=\, \sum\limits_{\,j \,\in\, J}\,v_{j}^{\,2}\,T^{\,\ast}\,P_{\,W_{j}}\,\Lambda_{j}^{\,\ast}\,\Lambda_{j}\,P_{\,W_{j}}\,U\,f.\]
By Cauchy-Schwartz inequality, for each \,$f \,\in\, H$, we have
\begin{align*}
\|\,f\,\|^{\,4} &\,=\, \left(\,\left<\,f,\, f\,\right>\,\right)^{\,2} \,=\, \left(\,\left<\,\sum\limits_{\,j \,\in\, J}\,v_{j}^{\,2}\,T^{\,\ast}\,P_{\,W_{j}}\,\Lambda_{j}^{\,\ast}\,\Lambda_{j}\,P_{\,W_{j}}\,U\,f,\, f\,\right>\,\right)^{\,2}\\
&=\, \left(\,\sum\limits_{\,j \,\in\, J}\,v_{j}^{\,2}\,\left<\,\Lambda_{j}\,P_{\,W_{j}}\,U\,f,\, \Lambda_{j}\,P_{\,W_{j}}\,T\,f\,\right>\,\right)^{\,2}\\
&\leq\, \sum\limits_{\,j \,\in\, J}\,v_{j}^{\,2}\,\left\|\,\Lambda_{j}\,P_{\,W_{j}}\,T\,f\,\right\|^{\,2}\,\sum\limits_{\,j \,\in\, J}\,v_{j}^{\,2}\,\left\|\,\Lambda_{j}\,P_{\,W_{j}}\,U\,f\,\right\|^{\,2} \\
& \,=\,\sum\limits_{\,j \,\in\, J}\,v_{j}^{\,2}\,\left<\,\Lambda_{j}\,P_{\,W_{j}}\,T\,f,\, \Lambda_{j}\,P_{\,W_{j}}\,T\,f\,\right> \,\sum\limits_{\,j \,\in\, J}\,v_{j}^{\,2}\,\left<\,\Lambda_{j}\,P_{\,W_{j}}\,U\,f,\, \Lambda_{j}\,P_{\,W_{j}}\,U\,f\,\right> \\
&\leq\, B\,\|\,f\,\|^{\,2}\,\sum\limits_{\,j \,\in\, J}\,v_{j}^{\,2}\,\left<\,\Lambda_{j}\,P_{\,W_{j}}\,U\,f,\, \Lambda_{j}\,P_{\,W_{j}}\,U\,f\,\right>\\
&\Rightarrow\,\dfrac{1}{B}\,\|\,f\,\|^{\,2} \,\leq\, \sum\limits_{\,j \,\in\, J}\,v_{j}^{\,2}\,\left<\,\Lambda_{j}\,P_{\,W_{j}}\,U\,f,\, \Lambda_{j}\,P_{\,W_{j}}\,U\,f\,\right>.  
\end{align*}
On the other hand, for each \,$f \,\in\, H$, we have
\begin{align*}
&\sum\limits_{\,j \,\in\, J}\,v_{j}^{\,2}\,\left<\,\Lambda_{j}\,P_{\,W_{j}}\,U\,f,\, \Lambda_{j}\,P_{\,W_{j}}\,U\,f\,\right>\\
&\hspace{.45cm} \,=\, \sum\limits_{\,j \,\in\, J}\,v_{j}^{\,2}\,\left<\,\Lambda_{j}\,P_{\,W_{j}}\,T\,T^{\,-\, 1}\,U\,f,\, \Lambda_{j}\,P_{\,W_{j}}\,T\,T^{\,-\, 1}\,U\,f\,\right>\\
&\hspace{.45cm}\leq\, B\,\left\|\,T^{\,-\, 1}\,U\,f\,\right\|^{\,2} \,\leq\, B\,\left\|\,T^{\,-\, 1}\,\right\|^{\,2}\,\|\,U\,\|^{\,2}\,\|\,f\,\|^{\,2}.
\end{align*}
Thus, \,$\Lambda_{T}$\, is a \,$(\,U,\,U\,)$-controlled $g$-fusion frame for  \,$H$.\,Similarly, it can be shown that if \,$\Lambda_{T}$\, is a \,$(\,U,\,U\,)$-controlled $g$-fusion Bessel sequence in \,$H$\, then \,$\Lambda_{T}$\, is also a \,$(\,T,\,T\,)$-controlled $g$-fusion frame for \,$H$.   
\end{proof}

\begin{definition}
Let \,$\Lambda_{T} \,=\, \left\{\,\left(\,W_{j},\, \Lambda_{j},\, w_{j}\,\right)\,\right\}_{j \,\in\, J}$\, and \,$\Gamma_{U} \,=\, \left\{\,\left(\,V_{j},\, \Gamma_{j},\, v_{j}\,\right)\,\right\}_{j \,\in\, J}$\, be \,$(\,T,\,T\,)$-controlled and \,$(\,U,\,U\,)$-controlled $g$-fusion Bessel sequences in \,$H$\, with bounds \,$D_{\,1}$\, and \,$D_{\,2}$, respectively.\,Then the operator \,$S_{T\,\Lambda\,\Gamma\,U} \,:\, H \,\to\, H$, defined by
\[S_{T\,\Lambda\,\Gamma\,U}\,(\,f\,) \,=\, \sum\limits_{\,j \,\in\, J}\,v_{j}\, w_{j}\,T^{\,\ast}\,P_{\,W_{j}}\,\Lambda^{\,\ast}_{j}\,\Gamma_{j}\,P_{\,V_{j}}\,U\,f\; \;\forall\; f \,\in\, H,\] is called the frame operator for the pair of controlled g-fusion Bessel sequences \,$\Lambda_{T}$\, and \,$\Gamma_{U}$.
\end{definition}

It is easy to verify that \,$\Lambda_{T}$\, and \,$\Gamma_{U}$\, are also two \,$g$-fusion Bessel sequences in \,$H$.\,So, the frame operator  
\[S_{\Lambda\,\Gamma}\,f \,=\, \sum\limits_{\,j \,\in\, J}\,v_{\,j}\,w_{\,j}\,P_{\,W_{j}}\,\Lambda_{j}^{\,\ast}\,\Gamma_{j}\,P_{\,V_{j}}\,f\; \;\forall\; f \,\in\, H\]
for this pair of \,$g$-fusion Bessel sequences is well-defined and bounded.\,Therefore, we can write \,$S_{T\,\Lambda\,\Gamma\,U} \,=\, T^{\,\ast}\,S_{\Lambda\,\Gamma}\,U$.\,Thus, \,$S_{T\,\Lambda\,\Gamma\,U}$\, is well-defined and bounded.\,Also, for each \,$f,\, g \,\in\, H$, we have
\begin{align*}
\left<\,S_{T\,\Lambda\,\Gamma\,U}\,f,\; g\,\right> &\,=\, \left<\,\sum\limits_{\,j \,\in\, J}\,v_{j}\, w_{j}\,T^{\,\ast}\,P_{\,W_{j}}\,\Lambda^{\,\ast}_{j}\,\Gamma_{j}\,P_{\,V_{j}}\,U\,f,\, g\,\right>\\
&=\, \sum\limits_{\,j \,\in\, J}\, v_{\,j}\, w_{\,j}\, \left<\,f,\, U^{\,\ast}\,P_{\,V_{j}}\,\Gamma^{\,\ast}_{j}\,\Lambda_{j}\, P_{\,W_{j}}\,T\,g\,\right>\\
&=\,\left<\,f,\, \sum\limits_{\,j \,\in\, J}\,w_{\,j}\,v_{\,j}\,U^{\,\ast}\,P_{\,V_{j}}\,\Gamma^{\,\ast}_{j}\,\Lambda_{j}\, P_{\,W_{j}}\,T\,g\,\right> \,=\, \left<\,f,\, S_{U\,\Gamma\,\Lambda\,T}\,g\,\right>
\end{align*} 
and hence \,$S^{\,\ast}_{T\,\Lambda\,\Gamma\,U} \,=\, S_{U\,\Gamma\,\Lambda\,T}$.\\

In the next theorem, we give a characterization of a controlled $g$-fusion frame with the help of frame operator. 

\begin{theorem}
Let \,$T \,\in\, \mathcal{G}\,\mathcal{B}\,(\,H\,)$\, and \,$\Lambda_{T}$\, be a \,$(\,T,\,T\,)$-controlled $g$-fusion Bessel sequence in \,$H$.\,Then \,$\Lambda_{T}$\, is a \,$(\,T,\,T\,)$-controlled $g$-fusion frame for \,$H$\, if and only if there exist an operator \,$U \,\in\, \mathcal{G}\,\mathcal{B}\,(\,H\,)$\, and a \,$(\,U,\,U\,)$-controlled $g$-fusion Bessel sequence \,$\Gamma_{U}$\, such that \,$S_{U\,\Gamma\,\Lambda\,T} \,\geq\, m\,I_{H}$, for some \,$m \,>\, 0$.    
\end{theorem}

\begin{proof}
Let \,$\Lambda_{T}$\, be a \,$(\,T,\,T\,)$-controlled $g$-fusion frame with bounds \,$A$\, and \,$B$.\,Now, we take \,$U \,=\, T$, and \,$V_{j} \,=\, W_{j}$, \,$\Lambda_{j} \,=\, \Gamma_{j}$\, and \,$v_{j} \,=\, w_{j}$, for all \,$j \,\in\, J$.\,Thus, for each \,$f \,\in\, H$, we have
\begin{align*}
&\left<\,S_{T\,\Lambda\,\Lambda\,T}\,f,\, f\,\right> \,=\, \left<\,\sum\limits_{\,j \,\in\, J}\,w^{\,2}_{\,j}\,T^{\,\ast}\,P_{\,W_{j}}\,\Lambda_{j}^{\,\ast}\,\Lambda_{j}\,P_{\,W_{j}}\,T\,f,\, f\,\right>\\
&\hspace{2.5cm}=\, \sum\limits_{\,j \,\in\, J}\,w^{\,2}_{\,j}\,\left<\,\Lambda_{j}\,P_{\,W_{j}}\,T\,f,\, \Lambda_{j}\,P_{\,W_{j}}\,T\,f\,\right> \,\geq\, A\,\|\,f\,\|^{\,2}.
\end{align*}

Conversely, suppose that there exist an operator \,$U \,\in\, \mathcal{G}\,\mathcal{B}\,(\,H\,)$\, and a \,$(\,U,\,U\,)$-controlled $g$-fusion Bessel sequence \,$\Gamma_{U}$\, with bound \,$D$\, such that 
\[m\,\|\,f\,\|^{\,2} \,\leq\, \left<\,S_{U\,\Gamma\,\Lambda\,T}\,f,\, f\,\right>\; \;\forall\; f \,\in\, H,\; \;\text{for some $m \,>\, 0$}.\]
Thus, for each \,$f \,\in\, H$, we have
\begin{align*}
m\,\|\,f\,\|^{\,2} &\,\leq\, \left<\,S_{U\,\Gamma\,\Lambda\,T}\,f,\, f\,\right> \,=\, \sum\limits_{\,j \,\in\, J}\,\left<\,v_{\,j}\,w_{\,j}\,U^{\,\ast}\,P_{\,V_{j}}\,\Gamma^{\,\ast}_{j}\,\Lambda_{j}\,P_{\,W_{j}}\,T\,f,\, f\,\right>\\
&=\, \sum\limits_{\,j \,\in\, J}\,\left<\,w_{\,j}\,\Lambda_{j}\,P_{\,W_{j}}\,T\,f,\, v_{\,j}\,\Gamma_{j}\,P_{\,V_{j}}\,U\,f\,\right>\\
&\,\leq\, \left(\,\sum\limits_{\,j \,\in\, J}\, v_{\,j}^{\,2}\; \left\|\,\Gamma_{j}\, P_{\,V_{j}}\,U\,f\,\right\|^{\,2}\,\right)^{1 \,/\, 2}\,\left(\,\sum\limits_{\,j \,\in\, J}\, w_{\,j}^{\,2}\;  \left\|\,\Lambda_{j}\, P_{\,W_{j}}\,T\,f\,\right\|^{\,2}\,\right)^{1 \,/\, 2}\\
&\leq\, \sqrt{D}\,\|\,f\,\|\,\left(\,\sum\limits_{\,j \,\in\, J}\, w_{\,j}^{\,2}\, \left\|\,\Lambda_{j}\,P_{\,W_{j}}\,T\,f\,\right\|^{\,2}\,\right)^{1 \,/\, 2}\\
&\Rightarrow\, \dfrac{m^{\,2}}{D}\,\|\,f\,\|^{\,2} \,\leq\, \sum\limits_{\,j \,\in\, J}\,w_{\,j}^{\,2}\,\left<\,\Lambda_{j}\,P_{\,W_{j}}\,T\,f,\, \Lambda_{j}\,P_{\,W_{j}}\,T\,f\,\right>.
\end{align*}
Thus, \,$\Lambda_{T}$\, is a \,$(\,T,\,T\,)$-controlled $g$-fusion frame for \,$H$.     
\end{proof}

\begin{theorem}
Let \,$S_{T\,\Lambda\,\Gamma\,U}$\, be the frame operator for a pair of \,$(\,T,\,T\,)$-controlled and \,$(\,U,\,U\,)$-controlled \,$g$-fusion Bessel sequences \,$\Lambda_{T}$\, and \,$\Gamma_{U}$\, with bounds \,$D_{\,1}$\, and \,$D_{\,2}$, respectively.\,Suppose \,$\lambda_{\,1} \,<\, 1,\; \lambda_{\,2} \,>\, \,-\,1$\, such that for each \,$f \,\in\, H$, 
\[\left\|\,f \,-\, S_{T\,\Lambda\,\Gamma\,U}\,f\,\right\| \,\leq\; \lambda_{\,1}\, \|\,f\,\| \,+\, \lambda_{\,2}\, \left\|\,S_{T\,\Lambda\,\Gamma\,U}\,f \,\right\|.\] Then \,$\Gamma_{U}$\, is a \,$(\,U,\,U\,)$-controlled $g$-fusion frame for \,$H$.   
\end{theorem}

\begin{proof}
For each \,$f \,\in\, H$, we have
\begin{align}
&\left(\,1 \,-\, \lambda_{\,1}\,\right)\,\|\,f\,\| \,\leq\, \left(\,1 \,+\, \lambda_{\,2}\,\right)\,\|\,S_{T\,\Lambda\,\Gamma\,U}\,f \,\|\nonumber\\
&\Rightarrow\, \left(\,\dfrac{1 \,-\, \lambda_{\,1}}{1 \,+\, \lambda_{\,2}}\,\right)\,\|\,f\,\| \,\leq\, \|\,S_{T\,\Lambda\,\Gamma\,U}\,f \,\| \,=\, \sup\limits_{\|\,g\,\| \,=\, 1}\,\left<\,S_{T\,\Lambda\,\Gamma\,U}\,f,\, g\,\right>\nonumber\\
&=\, \sup\limits_{\|\,g\,\| \,=\, 1}\,\left<\,\sum\limits_{\,j \,\in\, J}\,v_{j}\, w_{j}\,T^{\,\ast}\,P_{\,W_{j}}\,\Lambda^{\,\ast}_{j}\,\Gamma_{j}\,P_{\,V_{j}}\,U\,f,\, g\,\right>\nonumber\\
&\,\leq\, \sup\limits_{\|\,g\,\| \,=\, 1}\,\left(\,\sum\limits_{\,j \,\in\, J}\, v_{\,j}^{\,2}\; \left\|\,\Gamma_{j}\, P_{\,V_{j}}\,U\,f\,\right\|^{\,2}\,\right)^{1 \,/\, 2}\,\left(\,\sum\limits_{\,j \,\in\, J}\, w_{\,j}^{\,2}\;  \left\|\,\Lambda_{j}\, P_{\,W_{j}}\,T\,f\,\right\|^{\,2}\,\right)^{1 \,/\, 2}\nonumber\\
&\leq\, \sqrt{D_{\,1}}\,\left(\,\sum\limits_{\,j \,\in\, J}\, v_{\,j}^{\,2}\,\left\|\,\Gamma_{j}\, P_{\,V_{j}}\,U\,f\,\right\|^{\,2}\,\right)^{1 \,/\, 2}\nonumber\\
&\Rightarrow\, \dfrac{1}{D_{\,1}}\,\left(\,\dfrac{1 \,-\, \lambda_{\,1}}{1 \,+\, \lambda_{\,2}}\,\right)^{\,2}\,\|\,f\,\|^{\,2} \,\leq\, \sum\limits_{\,j \,\in\, J}\, v_{\,j}^{\,2}\,\left<\,\Gamma_{j}\, P_{\,V_{j}}\,U\,f,\, \Gamma_{j}\, P_{\,V_{j}}\,U\,f\,\right>\; \;\;\forall\, f \,\in\, H.\label{eq7}
\end{align} 
Thus, \,$\Gamma_{U}$\, is a \,$(\,U,\,U\,)$-controlled $g$-fusion frame for \,$H$.
\end{proof}

\begin{theorem}\label{th4}
Let \,$S_{T\,\Lambda\,\Gamma\,U}$\, be the frame operator for a pair of \,$(\,T,\,T\,)$-controlled and \,$(\,U,\,U\,)$-controlled \,$g$-fusion Bessel sequences \,$\Lambda_{T}$\, and \,$\Gamma_{U}$\, with bounds \,$D_{\,1}$\, and \,$D_{\,2}$, respectively.\,Assume \,$\lambda \,\in\, (\,0 \,,\, 1\,)$\; such that 
\[\left\|\,f \,-\, S_{T\,\Lambda\,\Gamma\,U}\,f\,\right\| \,\leq\, \lambda\, \|\,f\,\|\;  \;\;\forall\; f \,\in\, H.\]
Then \,$\Lambda_{T}$\, is a \,$(\,T,\,T\,)$-controlled $g$-fusion frame for \,$H$\, and \,$\Gamma_{U}$\, is a \,$(\,U,\,U\,)$-controlled $g$-fusion frame for \,$H$. 
\end{theorem}

\begin{proof}
By putting \;$\lambda_{\,1} \,=\, \lambda\; \;\text{and}\; \; \lambda_{\,2} \,=\, 0$\; in (\ref{eq7}), we get 
\[ \dfrac{(\,1 \,-\, \lambda\,)^{\,2}}{D_{\,1}}\, \|\,f\,\|^{\,2} \,\leq\, \sum\limits_{\,j \,\in\, J}\, v_{\,j}^{\,2}\,\left<\,\Gamma_{j}\, P_{\,V_{j}}\,U\,f,\, \Gamma_{j}\, P_{\,V_{j}}\,U\,f\,\right>\]and therefore \,$\Gamma_{U}$\, is a \,$(\,U,\,U\,)$-controlled $g$-fusion frame for \,$H$.\,Now, for each \,$f \,\in\, H$, we have 
\begin{align*}
&\left(\,1 \,-\, \lambda\,\right)\,\|\,f\,\| \,\leq\, \|\,S_{T\,\Lambda\,\Gamma\,U}^{\,\ast}\,f \,\| \,=\, \sup\limits_{\|\,g\,\| \,=\, 1}\,\left<\,S^{\,\ast}_{T\,\Lambda\,\Gamma\,U}\,f,\, g\,\right> \,=\, \sup\limits_{\|\,g\,\| \,=\, 1}\,\left<\, S_{U\,\Gamma\,\Lambda\,T}\,f,\, g\,\right>\\
&=\, \sup\limits_{\|\,g\,\| \,=\, 1}\,\left<\,\sum\limits_{\,j \,\in\, J}\,v_{j}\, w_{j}\,U^{\,\ast}\,P_{\,V_{j}}\,\Gamma^{\,\ast}_{j}\,\Lambda_{j}\,P_{\,W_{j}}\,T\,f,\, g\,\right>\\
&\,\leq\, \sup\limits_{\|\,g\,\| \,=\, 1}\,\left(\,\sum\limits_{\,j \,\in\, J}\, v_{\,j}^{\,2}\; \left\|\,\Gamma_{j}\, P_{\,V_{j}}\,U\,g\,\right\|^{\,2}\,\right)^{1 \,/\, 2}\,\left(\,\sum\limits_{\,j \,\in\, J}\, w_{\,j}^{\,2}\,\left\|\,\Lambda_{j}\, P_{\,W_{j}}\,T\,f\,\right\|^{\,2}\,\right)^{1 \,/\, 2}\\
&\leq\, \sqrt{D_{\,2}}\,\left(\,\sum\limits_{\,j \,\in\, J}\, w_{\,j}^{\,2}\,\left\|\,\Lambda_{j}\,P_{\,W_{j}}\,T\,f\,\right\|^{\,2}\,\right)^{1 \,/\, 2}\\
&\Rightarrow\, \sum\limits_{\,j \,\in\, J}\,w_{j}^{\,2}\,\left<\,\Lambda_{j}\,P_{\,W_{j}}\,T\,f,\, \Lambda_{j}\,P_{\,W_{j}}\,T\,f\,\right>  \,\geq\, \dfrac{(\,1 \,-\, \lambda\,)^{\,2}}{D_{\,2}}\; \|\,f\,\|^{\,2}\; \;\;\forall\, f \,\in\, H.
\end{align*}
Hence, \,$\Lambda_{T}$\, is a \,$(\,T,\,T\,)$-controlled $g$-fusion frame for \,$H$.\,This completes the proof.  
\end{proof}

\end{document}